\documentclass[11pt]{article}
\usepackage{amsmath,amsfonts,amsthm,amsbsy}
\usepackage{wrapfig}
 \usepackage[pdftex]{graphicx}
 \DeclareGraphicsExtensions{.pdf}
\numberwithin{equation}{section}
\newcommand{\set}[1]{ \{#1\} }

\newcommand{\pos}{\mathbb{P}}
\newcommand{\nat}{\mathbb{N}_0}
\newcommand{\integer}{\mathbb{Z}} 
\newcommand{\real}{\mathbb{R}}

\def\al{\alpha}
\newcommand{\be}{\beta}

\newcommand{\ka}{\kappa}

\def\tb{\mathbf{\tau}}
\newcommand{\om}{\omega}

\newcommand{\Om}{\Omega}
\newcommand{\bb}{\bar{\beta}}
\newcommand{\ab}{\bar{\alpha}}

\newcommand{\ra}{\rangle}
\newcommand{\la}{\langle}

\newcommand{\figref}[1]{Figure \ref{#1}}
\newcommand{\defref}[1]{Definition \ref{#1}}
\newcommand{\propref}[1]{Proposition \ref{#1}}
\newcommand{\lemref}[1]{Lemma \ref{#1}}
\newcommand{\thmref}[1]{Theorem \ref{#1}}

\newcommand{\tm}{\mathcal{P}}  % transition matrix
\newcommand{\ssv}{\vec{P}_S} % stationary state vector

\theoremstyle{plain}

\newtheorem{theorem}{Theorem}

\newtheorem{lemma}{Lemma}

\newtheorem{prop}{Proposition}

\newtheorem{definition}{Definition}
\theoremstyle{definition}
\newtheorem{remark}{Remark}

\newcommand{\stepSet}{\mathcal{S}}
\newcommand{\edgeseq}{\mathcal{E}}
\newcommand{\suchthat}{\, |\,} 
\newcommand{\st}{\, |\,} 
 
\begin{document} 
%
% Switches for changing form .eps files to/from .pdf files.
%
% \DeclareGraphicsExtensions{.eps}
% \DeclareGraphicsExtensions{.pdf}
%

%
%\pagenumbering{roman}

\title{Simple Asymmetric Exclusion Model and Lattice Paths: Bijections and Involutions}

\author{R. Brak\dag and J. W. Essam\ddag \\
		\vspace{1em}\\
        \thanks{email: {\tt r.brak@ms.unimelb.edu.au, J.Essam@rhul.ac.uk}}
         \dag Department of Mathematics,\\
         The University of Melbourne\\
         Parkville, \\
         Victoria 3052,\\
         Australia\\
\vspace{1em}\\
         \ddag Department of Mathematics,\\
         Royal Holloway College, University of London,\\
         Egham,\\
         Surrey TW20 0EX,\\
         England.}

\date{\today} 
 
\maketitle 
 
\newpage
\begin{abstract} 
We study the combinatorics of the change of basis of three  representations of the  stationary state algebra of the two parameter simple asymmetric exclusion process. Each of the  representations considered correspond to a different set of weighted lattice paths which, when summed over, give the stationary state probability distribution. We show that all three sets of paths are combinatorially related via  sequences of bijections and sign reversing involutions.
\end{abstract}

\vfill
\noindent{\bf Short title:} ASEP and Lattice Paths: Bijections and Involutions
 
 \bigskip


\newpage

% 

%\newpage
\section{Introduction}

\begin{figure}[ht]
\begin{center}
\includegraphics{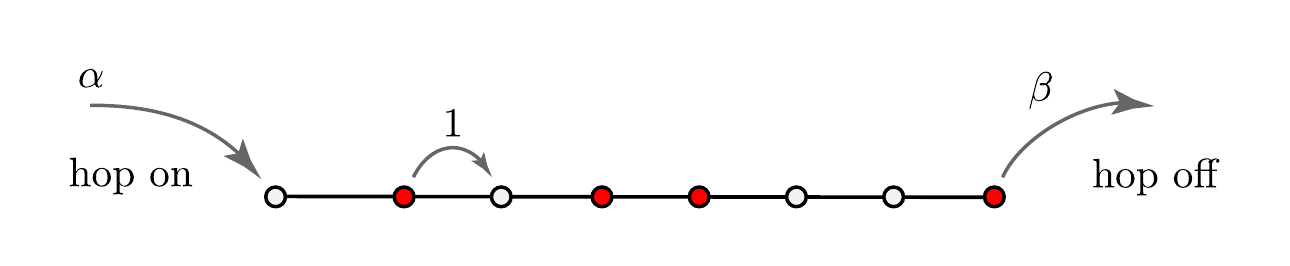} 
\caption{ASEP hopping model}\label{fig:hop}
\end{center}
\end{figure}

The Simple  Asymmetric Exclusion Process (ASEP) is a stochastic process defined by particles hopping along a line of length $L$ -- see \figref{fig:hop}.  Particles hop on to the line on the left with probability $\alpha$, off at the right with probability $\beta$ and between vertices to the right with unit probability with the constraint that only one particle can occupy a vertex. 
  The problem of readily computing the  stationary probability distribution was solved by Derrida et al \cite{derrida97}    with the introduction of the ``matrix product" Ansatz (see below) which provides an algebraic method of computing the stationary distribution. The ASEP 
and variations of it are a rich source of combinatorics: progress has been made in understanding the stationary distribution purely combinatorially \cite{Corteel05ve,duchi:2004ul, duchi:2004ep} and computing the stationary distribution has been shown to be equivalent to solving various lattice path problems \cite{Brak2004vf} or   permutation tableaux
\cite{williams:2006ub}. A recent review of the Asymmetric Exclusion Process may be found in  Blythe  and  Evans \cite{Blythe:2007uq}.

As explained in detail below, the matrix product Ansatz expresses the stationary distribution of a given state as a matrix product (the exact form of the product depends on the state). The matrices arise as representations of the DEHP  algebra. The paper by Derrida et al \cite{derrida97} originally found three different representations.  As shown by Brak and Essam \cite{Brak2004vf},   each matrix representation can be interpreted as a transfer  matrix (see \cite{stanley:1997vw} section 4.7) for a different lattice path model. Computing the stationary distribution is thus translated into finding certain lattice path weight polynomials. 
 
\begin{figure}[ht]
	\begin{center}
		\includegraphics{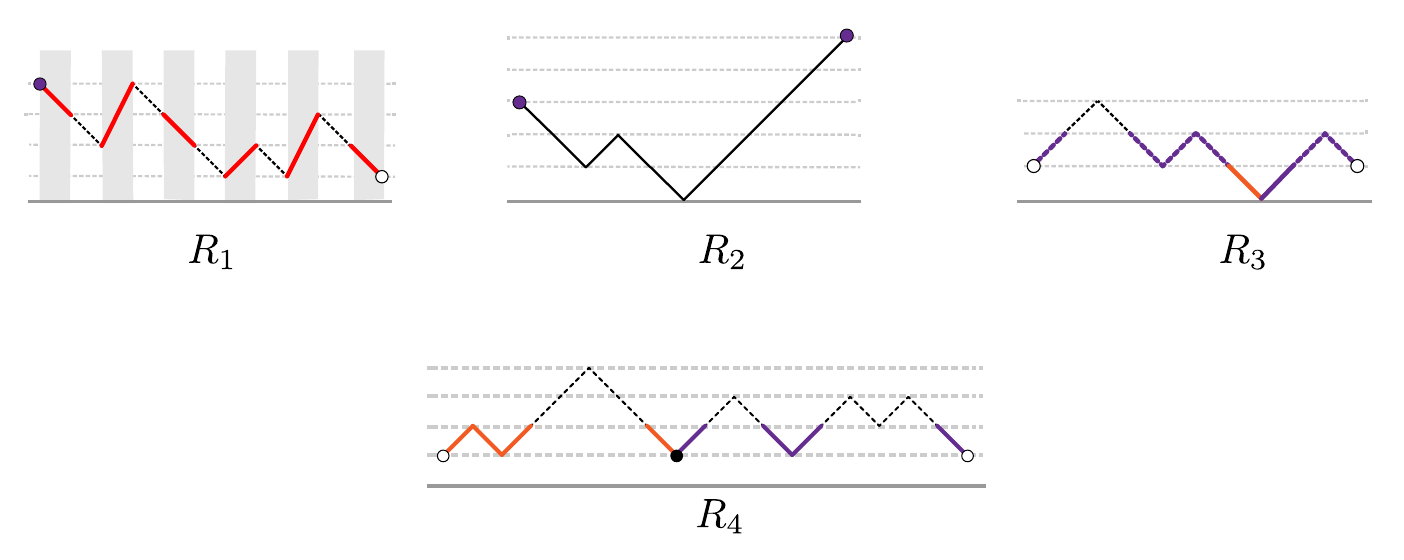}
		\caption{The lattice path models associated with the three algebra representations and the $R_4$, `canonical' representation path.}
		\label{fig:allbij}
	\end{center}
\end{figure}

Each of the three lattice path models are quite different (see - \figref{fig:allbij}) however \emph{they all have the same weight polynomials} (as they must since they all correspond to the same stationary probability). 
Our primary interest in this paper is to shown how this arises combinatorially.
This will be done by showing that all   three path models are related by weight preserving bijections and involutions. Rather than enunciate the three possible connections between the three paths  we rather show how they biject  to   a fourth ``canonical'' path model    -- see \figref{fig:allbij}. 

The primary consequences of theses connections are two-fold. Firstly the canonical path model provides a new  representation of the DEHP algebra and secondly, since each of these lattice paths arise from  representations of the DEHP algebra  the bijections between the different representations correspond algebraically to similarity transformations between the representations. Although we don't do so in this paper, it would be interesting to see how (if at all) the bijections are related to the similarity matrices themselves.

An additional  interest of the   canonical path model is that it can be interpreted as an interface polymer model. This polymer model has recently been used \cite{brak2004je} to gain a new understanding of how \emph{equilibrium} models in statistical mechanics are imbedded in \emph{non-equilibrium} process.
%
%===========================
\section{Markov chain and ASEP algebra}
%===========================
We now define the ASEP and briefly explain the Matrix product Anstaz.
The state  of the chain, $\tb=(\tau_1,\dots,\tau_L)\in (0,1)^L$, is determined by the particle occupancy
\begin{equation}
\tau_i=\begin{cases} 1&  \text{if vertex $i$ is occupied}\\ 
0 &\text{otherwise}
\end{cases}
\end{equation}
The transition matrix, $\tm$ has  elements,  
\begin{itemize}
\item  Hopping on: $ \tm\left( (0,\dots,\tau_L,),(1,\dots,\tau_L)\right)=\al$
\item  Hopping off: $ \tm\left( (\tau_1,\dots,\tau_{L-1},1 ),(\tau_1,\dots,\tau_{L-1},0)\right)=\beta$
\item  Right hopping: $ \tm\left( (\tau_1,\dots,\tau_{i}, \dots, \tau_L),(\tau_1,\dots,1-\tau_{i}, \dots, \tau_L)\right)=1$, for $\tau_i=1$, $1\le i<L$.
 
\end{itemize}
All other elements of $\tm$ are zero except the diagonals for which $\tm(\tau,\tau)=1-\sum_{\tau'\in (0,1)^L, \tau'\ne \tau} \tm(\tau',\tau)$
%
%  \begin{figure}[ht]
% \begin{center}
% \includegraphics[width=6cm]{\figPathSketch{sg_pasep.pdf}
% \caption{State space and transition probabilities for the $L=2$ chain.}
% \label{fig:state}
% \end{center}
% \end{figure}
% % 
The primary object we wish to determine is the stationary state vector $\ssv$ determined by\\ $\tm\ssv =0$. 
%The state diagram for $L=2$ is shown in \figref{fig:state}.  
 Derrida et al\cite{derrida97}, have shown that the stationary state vector
could be written as a  matrix product Ansatz, in particular they show the following.
%============= Theorem ========
\begin{theorem}\cite{derrida97} Let $D$ and $E$ be matrices then the components of the stationary state vector are given by 
\begin{equation}
P_{S}( \tau) =\frac{1}{Z_{L}}\,   W \, \left[\prod_{i=1}^{L}(\tau_{i}D+(1-\tau_{i})E)\right]\, V  \label{prob}
\end{equation}
with normalisation    $Z_{L}$ given by
\begin{equation} 
Z_{L} =  W  (D E )^{L} V  
\label{eq:norm}
\end{equation}
provided that $D$ and $E$  satisfy the DEHP algebra
\begin{subequations}\label{algebra}
\begin{align}
\label{eqn:de}
    D+E&=DE   
\intertext{and  $  W $ and $ V $ are the left and right eigenvectors}
      W E&=\frac{1}{\alpha}   W ,\qquad
     D  V  =\frac{1}{\beta } V .
\end{align}
\end{subequations}
\end{theorem}
%============= Theorem ========
These equations are sufficient to determine $P_{S}( \tau)$ algebraically.  Derrida et al \cite{derrida97} also gave several matrix representations of $D$ and $E$ and the vectors $|V \ra$ and $\la W |$,  any one of which may also used to determine $P_{S}( \tau)$. 

The three representations found by  Derrida et al \cite{derrida97}   are conveniently expressed in terms if the variables
\begin{subequations}\label{eq:vars}
\begin{align}
\ab&=1/\al\\
\bb&=1/\be\\
c&=\ab-1\\
d&=\bb-1\\
\ka^{2}&=\ab+\bb-\ab \bb=1-cd
\end{align}
\end{subequations}
and are as follows.

\paragraph{Representation I}
\begin{equation}
 		D_{1}= \begin{pmatrix}\bb&\bb&\bb&\bb&\bb&\cdots\\
		0&1&1&1&1&\cdots\\
		0&0&1&1&1&\cdots\\
		0&0&0&1&1&\cdots\\
		\vdots&\vdots&\vdots&\vdots&
		\end{pmatrix}  \qquad
		E_{1}= \begin{pmatrix}0&0&0&0&0&\cdots\\
		1&0&0&0&0&\cdots\\
		0&1&0&0&0&\cdots\\
		0&0&1&0&0&\cdots\\
		\vdots&\vdots&\vdots&\vdots&
		\end{pmatrix}  
		\label{eq:rep1}
\end{equation}
\begin{equation}
		 W_{1} =(1,\ab,\ab^{2},\ab^{3},\ldots)\qquad
		 V_{1}  = (1,0,0,0,\ldots)^{T}
		\label{eq:repv1}
\end{equation}
\paragraph{Representation II}
\begin{equation}
 		D_{2}= \begin{pmatrix}1&1&0&0&0&\cdots\\
		0&1&1&0&0&\cdots\\
		0&0&1&1&0&\cdots\\
		0&0&0&1&1&\cdots\\
		\vdots&\vdots&\vdots&\vdots&
		\end{pmatrix}  \qquad
		E_{2}= \begin{pmatrix}1&0&0&0&0&\cdots\\
		1&1&0&0&0&\cdots\\
		0&1&1&0&0&\cdots\\
		0&0&1&1&0&\cdots\\
		\vdots&\vdots&\vdots&\vdots&
		\end{pmatrix}  
		\label{eq:rep2}
\end{equation}
\begin{equation}
		  W_{2} =\ka(1,c,c^{2},c^{3},\ldots)\qquad
		 V_{2}  = \ka(1,d,d^{2},d^{3},\ldots)^{T}
		\label{eq:repv2}
\end{equation}
\paragraph{Representation III}
\begin{equation}
 		D_{3}= \begin{pmatrix}\bb&\ka&0&0&0&\cdots\\
		0&1&1&0&0&\cdots\\
		0&0&1&1&0&\cdots\\
		0&0&0&1&1&\cdots\\
		\vdots&\vdots&\vdots&\vdots&
		\end{pmatrix}  \qquad
		E_{3}= \begin{pmatrix}\ab&0&0&0&0&\cdots\\
		\ka&1&0&0&0&\cdots\\
		0&1&1&0&0&\cdots\\
		0&0&1&1&0&\cdots\\
		\vdots&\vdots&\vdots&\vdots&
		\end{pmatrix}  
		\label{eq:rep3}
\end{equation}
\begin{equation}
		 W_{3} =(1,0,0,0,\ldots)\qquad
		 V_{3}  = (1,0,0,0,\ldots)^{T}
		\label{eq:repv3}
\end{equation}
Each of these three matrices can be interpreted as the ``transfer matrix"  for a certain set of lattice paths. 

We will use the usual notation for the set of real numbers $\real$, integers $\integer$, non-negative integers $\nat$,  positive integers $\pos$, $[n]=\{i\in \pos \st 1\le i\le n\}$ and $n\dots m=\{i\in \integer \st n\le i\le m\}$.

Let  $G=(V,A)$ be a  pseudo-digraph (ie.\ directed graph with loops) with vertex set $V$ and arc set $A$. Associate arc weights $W_A: A\to \real$ and vertex weights $W_V: V\to \real$ with $G$. Denote the weighted  pseudo-digraph  by $G(W_A,W_V)$.  The \textbf{ transfer matrix}, $T(G)$ associated with the digraph $G(W_A,W_V)$  is the weighted adjacency matrix $T(G)$ with elements $T(G)_{i,j} =W_A(v_i,v_j)$ for all $(v_i,v_j)\in A$.  The important property of the transfer matrix for us is that it generates \textbf{weighted random walks} on $G$. A random walk of length $t\in\nat$ from vertex $u$ to vertex $v$ on $G$ is the arc sequence $r(u,v)=a_1a_2\dots a_t$ with $a_i=(u_i,v_i)\in A$ such that  $v_{i}=u_{i+1}$ for all $i\in  \{1,\dots, t-1\}$ with $u_1=u$ and $v_t=v$.  From the random walk we construct the $t$-step
\textbf{weight polynomial}, $Z^{(G)}_t(u,v)$ defined  by
\begin{equation}\label{eq:wepo}
Z^{(G)}_t(u,v)=W_V(u) \left[\sum_{r\in \Om^{(G)}_t(u,v) }\, \prod_{i=1}^t W_A(a_i(r))\right] W_V(v)
\end{equation}
where $\Om^{(G)}_t(u,v)$ is the set of all $t$ step random walks on $G$ from $u$ to $v$ and $a_i(r)$ is the arc $a_i$ in walk $r$. If there are no length $t$ random walks from $u$ to $v$  then $Z^{(G)}(u,v)=0$. Thus the walks pick  up the weight of the initial and final vertices as well as the weights of all the arcs they step across. The weight polynomial is simply related to the weighted adjacency matrix as given by the following classical lemma.
\begin{prop}\label{lem:tmx}
Let $G=(V,A)$  be a directed pseudo-graph  with weighted adjacency matrix, $T$, then the $t$ step weight polynomial, \eqref{eq:wepo}, is given by 
\begin{equation}
Z^{(G)}_t(u,v) = W_V(u)\,   \left(T^t\right)_{u,v} W_V(v).
\end{equation}
\end{prop}
It is conventional to  spread the random walk out in ``time" when it is then referred to as a lattice path.     

\begin{definition}[Lattice Path]\label{def:lpat}
	A length $t$ \textbf{lattice path},  $p$, on $\Xi$  is a sequence of vertices $v_{0}v_{1}\ldots v_{t}$, with $v_{i}\in \Xi$ and  $v_{i}-v_{i-1}\in\stepSet_{i}$ for all  $i\in [t]$, where $\stepSet_{i}$ is the $i^\text{th}$ \textbf{step set}  which contains the set of allowed  $i^\text{th}$ steps. 
	The set $\Xi$ is usually $\integer\times\integer$ or $\integer\times\nat$. The \textbf{height  of a vertex}, $v=(x,y)$ is the $y$ value. For a particular path, $p$, denote the corresponding sequence of steps by $\edgeseq(p)=e_{1}e_{2}\ldots e_{t}$ with $e_{i}=(v_{i-1},v_{i})$ for all  $i\in [t]$. 
	The \textbf{height of a step} is the height of its left vertex. 
	The step, $e_i$ is in an \textbf{even column} or is an \textbf{even step}  (respect.\ \textbf{odd column} or \textbf{odd step}) if $i$ is even (respt.\ odd).  We will associate a \textbf{vertex weight} $W: v_i\to\real$ with the initial, $i=0$, and final, $i=t$, vertices of the paths, as well as a \textbf{step weight} $W: e_i\to \real$ with each step, $i\in [t]$ of the path.  A $t=0$ length path is the single vertex $v_0\in\Xi$. Denote the length of a path $p$ by $|p|$. 
	A \textbf{subpath} of length $k$ of a lattice path, $p$ starting at $u$, is the path defined by a subsequence of adjacent vertices, $ v_{i} v_{i+1} \ldots v_{i+k-1}  v_{i+k}$, of the lattice path $p$ with $v_i=u$. If the first vertex and last vertex of the subpath has height $h$ and all other vertices of the subpath have height greater or equal to $h$, then the subpath is called \textbf{$\mathbf{h}$-elevated}. 
\end{definition}

Given a digraph $G$ we associate (somewhat arbitrarily) a lattice path. The weighted adjacency matrix determines the step sets as follows: $S_i=\{(1,i-j)\st \text{$(v_i,v_j)\in A(G)$ for all $v_j\in V(G)$}\}$.
Note, the step sets thus defined depend on the labelling of the vertices -- usually a labelling is chosen such that adjacent vertices, as far as possible, are labelled sequentially ie.\ $u$ and $v$ are labelled $v_i$ and $v_{i+1}$ if  $(v_i,v_{i+1})\in A(G)$. The vertex weights of the path are same as the vertex weights of the random walk, similarly then step weights of the path are the same as the corresponding arc weights of the random walk.

We can now consider  the three matrix representations, \eqref{eq:rep1}, \eqref{eq:rep2} and \eqref{eq:rep3} in the context of transfer matrices. For the normalisation, \eqref{eq:norm}  since only the product $D_iE_i$ occurs  the associated digraph $G_i$ is bipartite with, say vertex partition $V_{D_i}$ and $V_{E_i}$. Thus, $D_i$ represents part of the adjacency matrix for the weighted arcs from vertices in $V_{D_i}$ to  vertices are $V_{E_i}$    ie.\ the rows of $D_i$ are labelled by the vertices of $V_{D_i}$ and the columns of $D_i$ are labelled by the vertices of $V_{E_i}$. Similarly,  the weighted arcs  from $V_{E_i}$ to  $V_{D_i}$ are given  by $E_i$.
% -- this is illustrated schematically in \figref{fig:bip3}. 
Thus, labelling the vertices of the digraphs with positive integers gives  the  adjacency matrix, $T_i$.
\begin{equation}\label{eq:adtn}
(T_i)_{r,c}=\begin{cases}
(D_i)_{r,c} & \text{if $r$ is odd and $c$ is even}\\
(E_i)_{r,c} & \text{if $r$ is even and $c$ is odd}
\end{cases}
\end{equation}
where $r,c\in\pos$. Note, since the matrices $D_i$ and $E_i$ are infinite, so is the associated digraph. The   vertex weights   $W_{V}(k)$ of vertex $k$ in each of the vertex partitions $V_{E_i}$ and  $V_{D_i}$ are taken from the components of the corresponding eigenvectors,
\begin{subequations}\label{eq:vweg}
\begin{align}
 W_{V_{D_i}}(k)&= (W_i)_k\\
 W_{V_{E_i}}(k)&= (V_i)_k 
\end{align}
 \end{subequations}
 where $W_i$ and $V_i$, $i\in [3]$ are given by equations \eqref{eq:repv1},\eqref{eq:repv2} and \eqref{eq:repv3} respectively.
%
%  \begin{figure}[ht]
% \begin{center}
% \includegraphics[height=5cm]{\figPathSketch{bipartite3.pdf}
% \caption{The bipartite digraph for representation three.}\label{fig:bip3}
% \end{center}
% \end{figure}
% 
We now have the following relationship between random walks on digraphs (or equivalently lattice paths) and  the normalisation.
\begin{theorem}\cite{Brak2004vf}\label{thm:rw1}
Let $G_1$, $G_2$ and $G_3$ be directed  graphs with respective weighted adjacency matrices $T_1$, $T_2$ and $T_3$ defined by \eqref{eq:adtn} and vertex weights defined by \eqref{eq:vweg}.  The normalisation $Z_{L}$ defined in \eqref{eq:norm}  for the two-parameter ASEP is then given by the three expressions
\begin{subequations}\label{eq:dg7}
\begin{align}
Z_{L}&=\sum_{k\ge0}\, Z^{(G_1)}_{2L}(2k+1,1)\label{eq:dg7a}\\
Z_{L}&=\sum_{k\ge0}\sum_{\ell\ge0}\,Z^{(G_2)}_{2L}(2k+1,2\ell+1)\label{eq:dg7b}\\
Z_{L}&=Z^{(G_3)}_{2L}(1,1)\label{eq:dg7c}
\end{align}
\end{subequations}
where $Z^{(G)}_t(u,v)$ is given by \lemref{lem:tmx} 
\end{theorem}

%!======================================================
%:Path definitions
\subsection{The Three Lattice Path Models}
%!======================================================
 
\begin{figure}[ht]
\begin{center}
\includegraphics[width=30em]{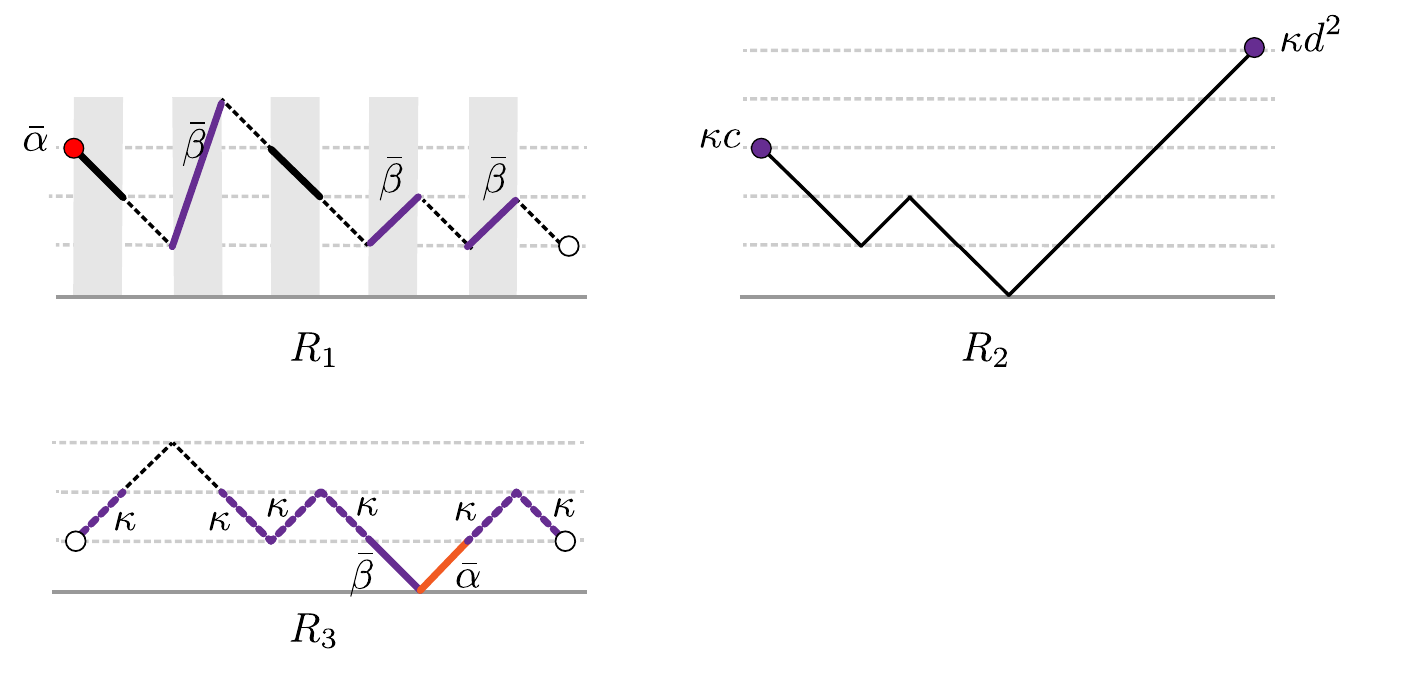}
\caption{\it An example of the three types of weighted paths, $R_1$, $R_2$  and $R_3$.}\label{fig:sw}
\label{fig:state}
\end{center}
\end{figure}

Associated with random walks on each of the three digraphs are lattice paths problems. Most of the lattice paths are  similar to Dyck paths. A \textbf{Dyck path} is a lattice path with step sets $S_i=\{(1,-1),(1,1)\}$ such that the height of the first vertex is the same as the height of the last vertex, and the height of all the remaining vertices is greater or equal to the the height of the first vertex. Examples of the first three types of lattice paths defined below are shown in \figref{fig:state}.
%
 %! R_1 =========  
\begin{definition}[$R_1$ paths]\label{def:r1p}
	$R_1$ paths are lattice paths on $\Xi=\integer\times\nat$ with step  sets
	\begin{equation}\label{eq:jspst1}
		\stepSet_{i}(R_1)=\left\{
		\begin{array}{ ll}
		   \{(1,-1)\}   & \text{for $i$ even (an `even  down step')}  \\
		    \{(1,2{k}-1) \st {k}\in \nat \}   &    \text{for $i$ odd (an `odd  (jump) step').}
		\end{array}
		\right.
	\end{equation}	
	with  $v_{0}(k)=(0,2k+1)$ for some $k\in \nat$ and $v_{2L}=(2L,1)$. Steps in $\{(1,2{k'}+1) \st {k'}\in \nat \}$   are called \textbf{jump up steps} and the \textbf{jump height}  is $2{k'}+1$. The $(1,-1)$ steps are called \textbf{odd down steps} (if $i$ is odd) or \textbf{even down steps} (if $i$ is even).   The weights associated with $R_1$ paths are 
	\begin{subequations}
	\begin{align}
		 W^{(1)}(v_0(k))&=\ab^{k}\\
		 W^{(1)}(v_{2L})&=1 \\
		 W^{(1)}(e_{i}) &=
			\begin{cases}
			 	\bb & \text{if  $e_{i}=\bigl((i-1,1),(i,2k'+2)\bigr)$, $k'\in \nat $ and  $i$   odd}\\
				1 & \text{otherwise}
			  \end{cases} 
	 \end{align}
	 \end{subequations}
 \end{definition}
Thus   $R_1$ paths start at some  odd height $y=2k+1$, every even step must be a down step, whilst an odd step may be a down step or a  step up an arbitrary (odd) jump height. The path must end at $({2L},1)$. 
Although the $R_1$ paths have a step from height one to height zero, there is no step  from height zero to one which combined with the constraint that the last step ends at height one means  $R_1$ paths have no vertices with height zero.
An example is shown in \figref{fig:sw}.
%
%! R_2 ============
%
\begin{definition}[$R_2$ paths]\label{def:r2p}
$R_2$ paths are lattice paths on $\Xi=\integer\times\nat$ with step  sets
\begin{equation} \label{eq:jspst2}
\stepSet_{i}(R_2)=   \{(1,-1), (1,1)\}      
 \end{equation}	
with  $v_{0}(k)=(0,2k+1)$ for some $k\in \nat$ and  $v_{{2L}}({k'})=(0,2{k'}+1)$ for some ${k'}\in \nat$.  The weights associated with $R_2$ paths are 
\begin{subequations}\label{eqs:w2}
\begin{align}
W^{(2)}(v_0(k))&=\ka\, c^{k}\\
W^{(2)}(v_{2L}({k'}))&=\ka\, d^{{k'}}\\
W^{(2)}(e_{i})& =1 \qquad\text{for all $i\in [t]$}
\end{align}
\end{subequations}
 \end{definition}
 Thus, $R_2$ paths are similar to  Dyck paths which start at height $2k+1$ and end at height $2{k'}+1$ with weights on the initial and final vertices. They  are also sometimes called ``rigged Ballot" paths. An example is shown in \figref{fig:sw}.
% 
%! R_3 ============
%
\begin{definition}[$R_3$ paths]\label{def:r32p}
$R_3$ paths are lattice paths on $\Xi=\integer\times\nat$ with step  set
\begin{equation} \label{eq:jspst3}
\stepSet_{i}(R_3)=   \{(1,-1), (1,1)\}      
 \end{equation}	
with  initial vertex $v_{0}=(0,1)$   and final vertex $v_{2L}=(2L,1)$. The weights associated with $R_3$ paths are 
\begin{subequations}
\begin{align}
 W^{(3)}(v_0)&=1,\\
 W^{(3)}(v_{2L})&=1,\\
  W^{(3)}(e_{i}) &=
			\begin{cases}
					\ka & \text{if $e_{i}=\bigl((i-1,1),(i,2)\bigr)$
					or $e_{i}=\bigl((i-1,2),(i,1)\bigr)$}\\
			   \bb & \text{if 
			   $e_{i}=\bigl((i-1,1),(i,0)\bigr)$ }\\
				 \ab & \text{if 
			   $e_{i}=\bigl((i-1,0),(i,1)\bigr)$ }\\
				 1 & \text{otherwise}
		  \end{cases}
 \end{align}
\end{subequations}
 \end{definition}
 Thus, $R_3$ are also similar to Dyck paths which start at height one and end at height one with weights on the first and second `levels'.  An example is shown in \figref{fig:sw}.  
 
\begin{figure}[ht]
\begin{center}
\includegraphics{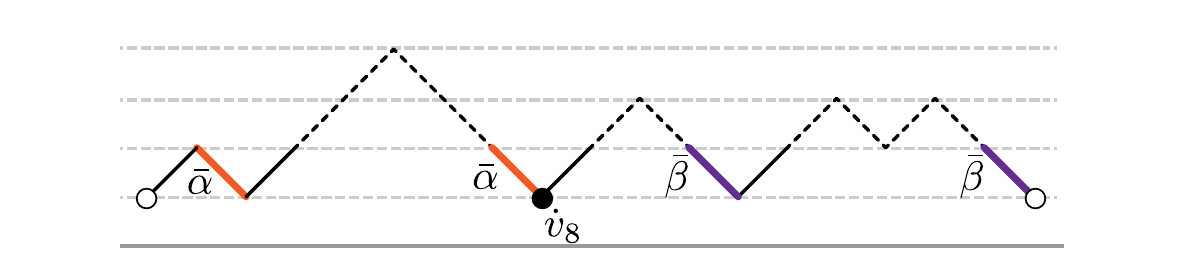} 
\caption{An example of an $R_4 $ path.}
\label{fig:r4p}
\end{center}
\end{figure}

 We now consider a fourth type of lattice path, which we will call $R_4 $ or `canonical' paths. They have   also been called \emph{one transit paths} \cite{brak2004je}  where they were used to model the behaviour of a polymer adsorbing on to an interface.
 % 
 % R4
 %
\begin{definition}[$R_4 $ paths]\label{def:r4p}
$R_4 $ paths are lattice paths on $\Xi=\integer\times\nat$ with step  sets
\begin{equation}  
	\stepSet_{i}(R_4 )=   \{(1,-1), (1,1)\}      
\end{equation}	
with  $v_{0}=(0,0)$, $v_{2L}=(2L,1)$ and one of the height one vertices marked. All vertices have height greater than zero. Denote the marked vertex with a dot, $\dot{v}$.  If $p$ is an $R_4 $ path and $p=v_0\dots\dot{v}_k\dots v_{2L}$, then the  weights associated with $R_4 $ paths are   
\begin{subequations}
	\begin{align}
	W^{(4)}(v_0)&=1,\\  
	W^{(4)}(v_{2L})&=1\\
	W^{(4)}(e_{i}) &=
		\begin{cases}
		   \ab & \text{if $e_{i}=\bigl((i-1,2),(i,1)\bigr)$ and $i\le k$}\\
		   \bb & \text{if  $e_{i}=\bigl((i-1,1),(i,2)\bigr)$ and $i>k$ }\\
		   1 & \text{otherwise}
		\end{cases}
	\end{align}
\end{subequations}
\end{definition} 
Note, $R_4$ paths are one-elevated, but there is a trivial bijection to zero-elevated paths, the one-elevation is  merely for convenience since most of the bijections and involutions discussed in this paper result naturally with the one-elevated form. Thus $R_4 $ paths are Dyck paths   with different weights to the left and right of the marked vertex.  An example is shown in \figref{fig:r4p}.  

The primary purpose of this paper is to provide a combinatorial proof of the following theorem.
\begin{theorem} \label{thm:ptg4}  Let $R_i(L)$ be the set  of $R_i$ paths of length $2L$.  
	The normalisation $Z_{L}$ defined in \eqref{eq:norm}  for the two parameter ASEP is   given by the four expressions
	\begin{subequations}\label{eq:lpf}
	\begin{align}
	 Z_{L}&=\sum_{p\in R_1(L)}\, W^{(1)}(p)\label{eq:lpf1}\\
	 Z_{L}&=\sum_{p\in R_2(L)}\, W^{(2)}(p)\label{eq:lpf2}\\
	 Z_{L}&=\sum_{p\in R_3(L)}\, W^{(3)}(p)\label{eq:lpf3}\\
	 Z_{L}&=\sum_{p\in R_4(L)}\, W^{(4)}(p)\label{eq:lpf4}
	\end{align}
	\end{subequations}
	where $W^{(i)}(p)$ is the weight of the path $p\in R_i(L)$. 
\end{theorem}
% 
%!============= Remark ========
\begin{remark}\label{rem:bhr}
Equations, \eqref{eq:lpf1}, \eqref{eq:lpf2} and \eqref{eq:lpf3} are essentially  those of   \thmref{thm:rw1} but stated in lattice paths form. Equation \eqref{eq:lpf4} is a new result.
\end{remark}
As an example, for $L=2$ the four expressions obtained are
\begin{align}\label{eq:hhn}
	Z_2(R_1)= &  \ab+\bb+ \ab \bb+\ab^2+\bb^2\\
	Z_2(R_2)= & 5+ 4\ka^2\frac{c}{1-cd}+4\ka^2\frac{d}{1-cd}+ \ka^2\frac{cd}{1-cd}+\ka^2\frac{c^2}{1-cd}+\ka^2\frac{d^2}{1-cd}\notag\\
	&= 5+4c+4d+cd+c^2+d^2\\
	Z_2(R_3)= & \ka^2+2\ab\bb\ka^2 +\ab^2\bb^2+\ka^4\\
	Z_2(R_4 )= & \ab+\bb+ \ab \bb+\ab^2+\bb^2.  \label{eq:exw4}
\end{align}
We make the following remarks based on the above example.
\begin{remark}\label{rem:wer4}
\begin{enumerate}

	\item In  equation  \eqref{eq:lpf1}  the final vertex $v_{2L}=(2L,1)$ is fixed, thus, even though arbitrary high jumps steps are permitted the constraint of having to return to height one results in a finite number of contributing paths,  hence the number of configurations is finite for fixed $L$, thus \eqref{eq:lpf1} is a polynomial in $\ab$ and $\bb$.
	
	\item  Equation \eqref{eq:lpf2}  is an infinite sum in  $c$ and $d$, but, as  will be shown, the infinite sum is always a simple geometric series giving rise to the  $(1-cd)^{-1}$ factor which cancels the common factor of $\ka^2=1-cd$  resulting in a polynomial in $c$ and $d$.   
	
	\item   Equation  \eqref{eq:lpf3}  is a polynomial in $\ab$, $\bb$ and $\ka^2$ arising from the three weights of the $R_3$ paths.
	
	\item  Equation  \eqref{eq:lpf4}  is a polynomial in $\ab$, $\bb$ -- the same polynomial as \eqref{eq:lpf1}, but arises from a completely different set of paths. 
	
\end{enumerate}
\end{remark}
The combinatorial proof shows how the four polynomials are connected  and how the $1-cd$ factor arises in the $R_2$  expression -- \eqref{eq:lpf2}. Combinatorially they are related by involutions and/or bijections to  the fourth  \eqref{eq:lpf4}. 
Thus there are three major parts to the proof each shows the connections between the three sets of paths and the $R_4 $ paths.
% ie.\  $R_1\to R_4 $, $R_2\to R_4 $ $R_3\to R_4 $. 
%There are several steps in each of the three cases which are summarised in \figref{fig:pst}
% %
%  \begin{figure}[ht]
% \begin{center}
% %\includegraphics{\figPathSketch{UoutlineTest.pdf}
% \includegraphics[width=30em]{\figPathSketch{bijecSteps.pdf}
% \caption{A summary of the steps used in the various combinatorial proofs}
% \label{fig:pst}
% \end{center}
% \end{figure}
% 

In each of the the  proofs  we will need to go via several different types of path before getting to $R_4 $ paths. 
Any type  of path  between an $R_i$ path   and the $R_4 $ path  will be labelled $R_{i}^{j}$, being the $j^\text{th}$ paths on route from $R_i$ paths to $R_4 $ paths.

Most of the proofs are  constructed by factoring the paths into certain subpaths.   We anticipate this by factoring $R_4 $ paths into Dyck subpaths by representing the lattice path using an alphabet.  
%===========
\begin{figure}
	\begin{center}
		\includegraphics[width=30em]{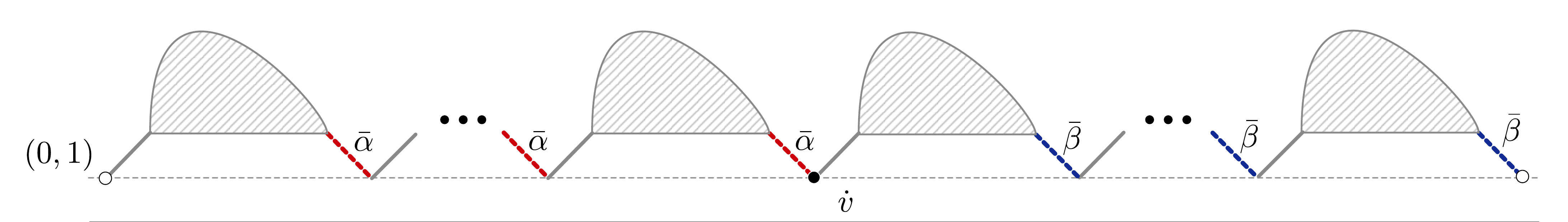}
		\caption{Schematic representation of the $D$-factorisation of a $R_4$ path also showing the $\ab$ weights to the left of all the $\bb$ weights (which defines the marked vertex $\dot{v}$).}
		\label{fig:Dfac}
	\end{center}
\end{figure}
%======= FIG ===========
Denote an up step by $u$ and a down step by $d$. If we scan the word representing the path from left to right noting the step associated with each time the path returns to height one (ie.\ a down step from height two to height one) then we have the following classical factorisation proposition  (illustrated schematically in \figref{fig:Dfac}.
).
\begin{prop}\label{lem:catfac}
Let $p\in R_4$ with weight $\ab^k\bb^{k'}$, then $p$ can be written in the   form
\begin{equation}
	p=\left[\prod_{i=1}^k u\, D_i\, d\right]\left[ \prod_{j=1}^{k'} u\, D_{j+k'}\, d\right]
\end{equation}
where $D_i$ is $2$-elevated (see \defref{def:lpat}) Dyck path. The weight  of a $d$ step  in the first factor is $\ab$ and the weight  of a $d$ step  in the second factor is $\bb$. If either $k=0$ or  ${k'}=0$   then corresponding product is absent. 
\end{prop}
We will refer to the above factorised form as  the \textbf{$D$-factorisation}.

%! |||||||||||||||||||||||||||||||||||||||||||||

%! R_1 -> R_4
\subsection{Proof of Equivalence of the $R_1$ and $R_4 $ path representations}
%! ***************************************************************

\newcommand{\auxI}[1]{A^{(#1)}_1} % Auxiliary path type 
We need to show
\begin{equation}\label{eq:ppf6}
 \sum_{p\in R_1}\, W^{(1)}(p) =\sum_{p\in R_4}\, W^{(4)}(p)\\
\end{equation}
where the paths in the sum are all of length $2L$.

The proof is by  bijection and proceeds in two stages, the first stage uses an  elevated subpath factorisation to biject to $R_1^1$ paths (defined below) and the second stage bijects the $R_1^1$ paths  to the $R_4$ paths.

\newcommand{\sed}{d_e\,} % even down 
\newcommand{\sju}[1]{u_{#1}\,}	% odd jump up
\newcommand{\sjd}{\bar{u}\,}	% odd jump down
\newcommand{\sdu}{u\,}	% Dyck up
\newcommand{\sdd}{d\,}	% Dyck down

\begin{figure}[ht]
	\begin{center}
	\includegraphics[width=11cm]{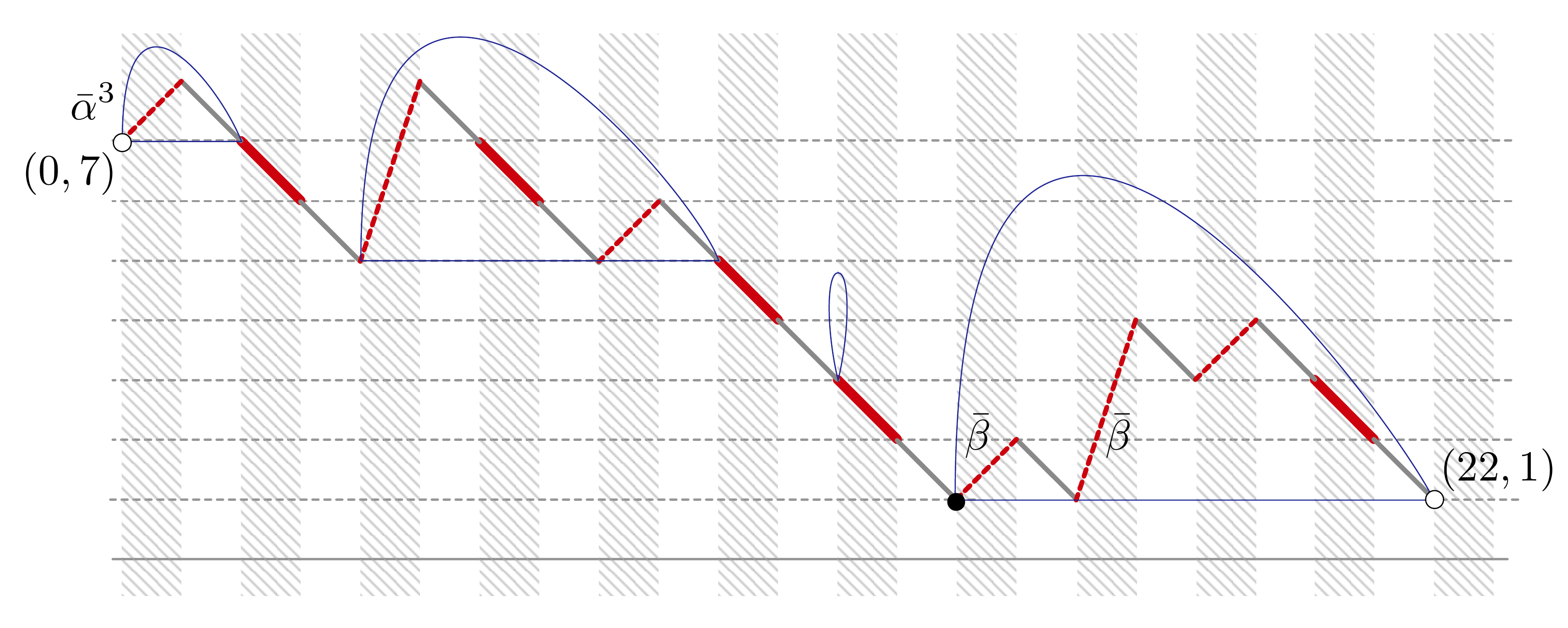}
	\caption{\it  The 22 step  $R_1$ path with  $\ab$ and $\bb$ weights given in equation \eqref{eq:exr1p} showing one level of $J$-factorisation.}
	\label{fig:r1ex}
	\end{center}
\end{figure}

When a path is represented by a step sequence (or word) the height of the initial vertex is not specified. Thus, if necessary we add the extra information by  representing the path as a pair  $k:w$ where $k$ is the height of the first vertex and $w$ is a word (or step sequence) in the alphabet $\{ \sjd, \sed,  \sju{2k+1}\}$ where $\sjd$ is a jump down step, $\sed$ an (even) down step and  $\sju{2k+1}$ a $2k+1$, $k\ge0$ jump step. 
  As an example, the $R_1$ path illustrated in \figref{fig:r1ex} is represented by
\begin{equation}\label{eq:exr1p}
7:  \sju{1} \sed \sjd \sed     \sju{3}\sed\sjd \sed    \sju{1}\sed\sjd \sed \sjd \sed     \sju{1}\sed       \sju{3} \sed \sju{1} \sed \sjd \sed \,. 
\end{equation}
%======= FIG ===========
%======= FIG ===========
We begin with a recursive factorisation  of the word representing a   path. The  recursion is simplest to state if the path  starts with  the   steps $\sed\sju{2k+1}\sed $ and ends with $\sjd$,  thus we define the factorisation  of    paths   in the set 
\begin{equation}
R'_1=\set{ 2: \sed\sju{2k+1}\sed \,\cdot w \cdot \sjd\st 2k+1 : w \in R_1,\, k\ge 0 }
\end{equation}
from which we obtain the factorisation of the paths in $R_1$.

\newcommand{\rec}[1]{{(#1)}}

%======= FIG ===========
 \begin{figure}[ht]
	\begin{center}
	\includegraphics[width=11cm]{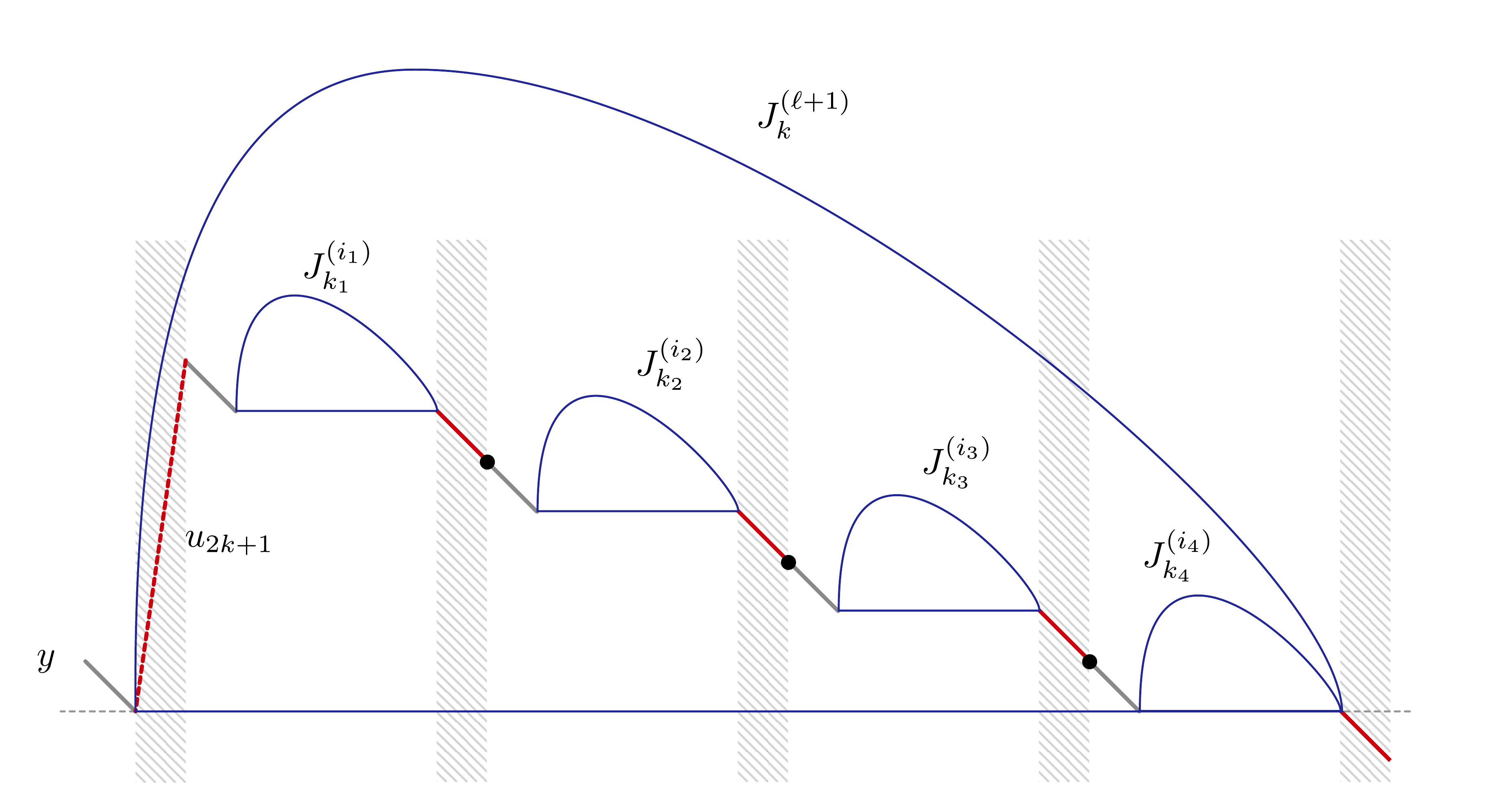}
	\caption{\it  A schematic representation of the $J$-factorisation ($k=3$ in this case) showing  one change in the level of recursion. }
	\label{fig:Jfac}
	\end{center}
\end{figure}
\noindent\emph{Stage 1:} The $R'_1$ paths   are factorised by reading the word from left to right starting after the initial $y:\sed\sju{2k+1}\sed$ prefix: the first time the path steps below the   height $y+2k-2j+2$ defines the `end' of the $J_{k_j}$ factor. This gives the following proposition.
%
%============= Prop ==============================
%
\begin{prop}
\label{lem:jfac}
Let $2k+1:w\in R_1$  , then   $ 2:w'\in R_1'$, $w'=\sed\sju{2k+1}\sed \,\cdot w \cdot \sjd$  has the recursive factorisation,
\begin{align} 
	w'=\sed J_k^\rec{\ell+1}\, \sjd =&
	\begin{cases}
		 \sed \sjd  										&\text{if $J_k^\rec{\ell+1}=\phi$}\\
		\displaystyle \sed \sju{2k+1} \prod_{j=1}^{k+1} \left[ \sed J_{k_j}^\rec{i_j} \sjd \right] 	& \text{if $J_{k}^\rec{\ell+1}  \ne\phi$}\label{eq:jfac2}
	\end{cases}\\
%\end{equation}
\intertext{where}
%\begin{equation}
 \ell = &\begin{cases} 
  			-1 & \text{if $J_k^\rec{\ell+1}=\phi$}\\
    		\max\set{i_j \st j=1\dots k+1}  & \text{otherwise}
		\end{cases}\label{eq:eq_lamx}
\end{align}
and $J_{m}^\rec{n} $ is a sub-path whose first step is a $2{m} +1$ jump step  and all vertices  of $J_{m}^\rec{n} $ are no lower than the first vertex  of the initial jump step of  $J_{m}^\rec{n} $. If $J_{m}^\rec{n} $ is empty  then the factor is denoted $J^{(0)}_\phi$. The superscript  $n$  on $J_{m}^\rec{n}$ denotes the level of recursion. The initial level is $J^{(0)}_\phi$.
\end{prop}
%============= Prop =================================
%
We will refer to the above form as  the \textbf{$J$-factorization}. The level of recursion is used primarily for induction proofs used below. One level change  of factorisation is illustrated schematically in \figref{fig:Jfac}.

\begin{figure}[ht]
	\begin{center}
	\includegraphics[width=30em]{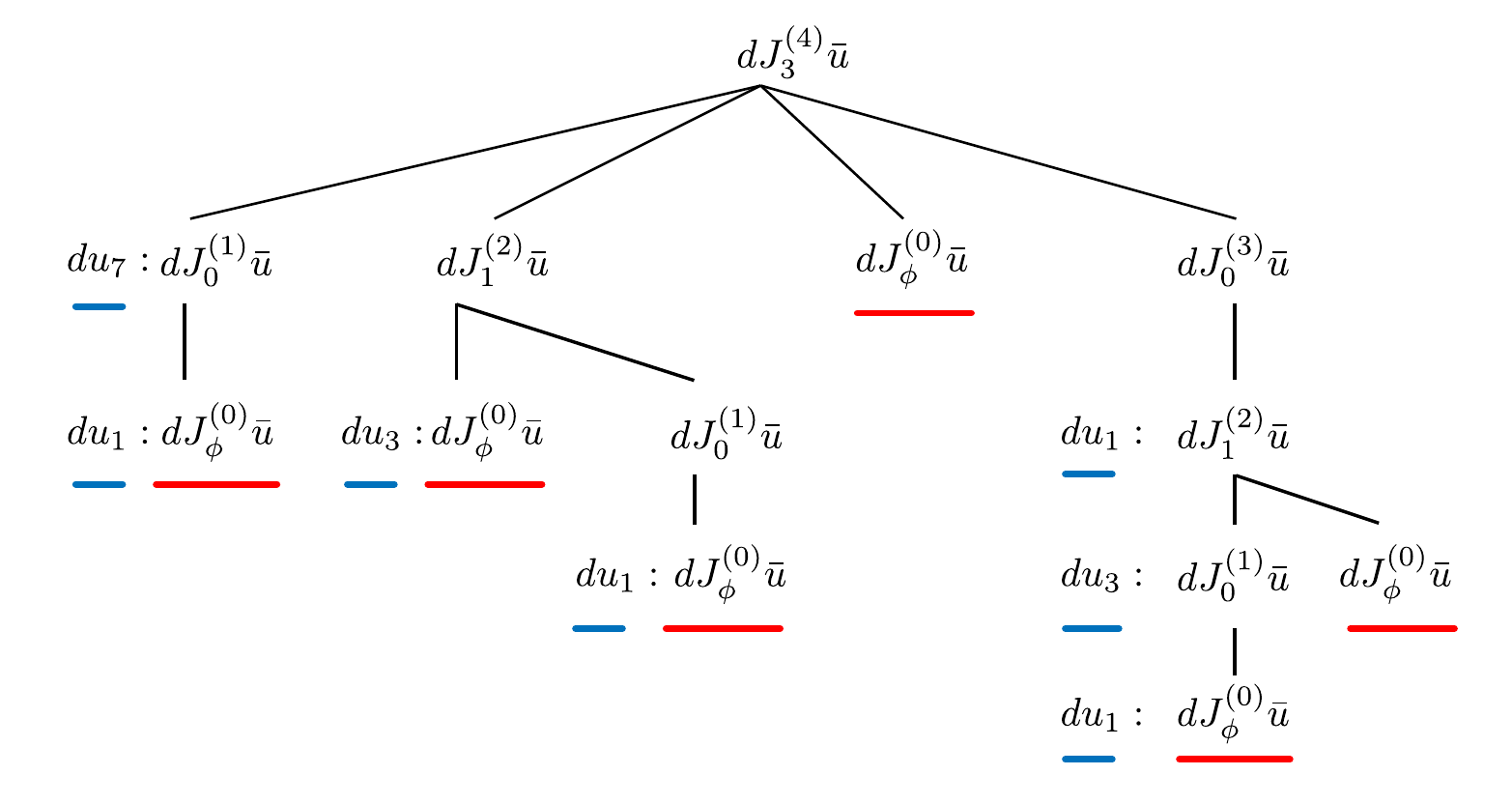}  
	\caption{\it The ordered planar tree representation of the $J$ factorisation of the path \eqref{eq:expp3} ($J^{(0)}_\phi$ is the empty factor). To read of the path the tree is traversed    `depth-first-search with left-to-right priority'  concatenating only the `$d u_{2k+1} $' (underlined blue)  prefixes and the `$d J \bar{u}$' factors (underlined red) on the leaves. }
	\label{fig:rroeeta1}
	\end{center}
\end{figure}

For example, the  path $p$, in \eqref{eq:exr1p} has the $J$-factorisation determined by factoring 
$\sed\sju{7} \sed\cdot p \cdot\sjd$ using \eqref{eq:jfac2}, as follows (the ``$\cdot$'' are only used to clarify the factorisation):
\begin{align}
\sed J_3\, \sjd &   =	 \sed\sju{7}\sed\cdot \sju{1} \sed \sjd \sed     \sju{3}\sed\sjd \sed    \sju{1}\sed\sjd \sed \sjd \sed     \sju{1}\sed \sjd    
		 \sju{3} \sed \sju{1} \sed \sjd \sed\cdot  \sjd\notag\\			
%		&=  \sed  \sju{7} \cdot  \bigl[ \sed J_0\, \sjd \bigr] \cdot \bigl[  \sed J_1\, \sjd \bigr]\cdot\bigl[\sed  \sjd\bigr] \circ \bigl[\sed J_0\, \sjd \bigr]  \notag \\	
%		&=  \sed  \sju{7} 
%		 \cdot 
%		 \bigl[ \sed \sju{1}\sed \, \sjd  \bigr] \cdot 
%		\bigl[  \sed  \sju{3}  \cdot 
%					\bigl[ \sed  \sjd  \bigr] \cdot 
%					\bigl[  \sed   J_0  \,   \sjd \bigr]\bigr] \cdot 
%		\bigl[ \sed  \sjd \bigr] \circ\notag \\
%		& \qquad 
%		\bigl[ \sed  \sju{1}  \cdot 
%					\bigl[ \sed J_1\, \sjd \bigr]  \bigr] \notag\\
%		&=  \sed  \sju{7} 
%		 \cdot 
%		 \bigl[ \sed \sju{1}\cdot [\sed \, \sjd]  \bigr] \cdot 
%		\bigl[  \sed  \sju{3}  \cdot 
%					\bigl[ \sed  \sjd  \bigr] \cdot 
%					\bigl[  \sed   J_0  \,   \sjd \bigr]\bigr] \cdot 
%		\bigl[ \sed  \sjd \bigr] \circ \notag\\
%		& \qquad 
%		\bigl[ \sed  \sju{1}  \cdot 
%					\bigl[ \sed  \sju{3} \cdot
%								\bigl[ \sed J_0\,  \sjd \bigr]  \cdot \bigl[\sed      \sjd 
%					\bigr]  
%		\bigr] \notag\\
		&=  \sed  \sju{7} 
		 \cdot 
		 \bigl[ \sed \sju{1}\sed \, \sjd  \bigr] \cdot 
		\bigl[  \sed  \sju{3}  \cdot 
					\bigl[ \sed  \sjd  \bigr] \cdot 
					\bigl[  \sed    \sju{1} \cdot [ \sed   \,   \sjd ] \bigr]\bigr] \cdot 
		\bigl[ \sed  \sjd \bigr] \cdot \notag\\
		& \qquad  \bigl[ \sed  \sju{1}  \cdot 
					\bigl[ \sed  \sju{3} \cdot
								\bigl[ \sed \sju{1} \cdot 
									[\sed\,  \sjd
									] 
								\bigr]
					\bigr]  \cdot 
					\bigl[
						\sed      \sjd 
					\bigr]  
			\bigr] \label{eq:expp3}
\intertext{thus, removing the prefix $\sed\sju{7} \sed $ and the $\sjd $ suffix, gives the factorisation of $p$ as,} 
p&= \bigl[   \sju{1}\sed \, \sjd  \bigr] \cdot 
		\bigl[  \sed  \sju{3}  \cdot 
					\bigl[ \sed  \sjd  \bigr] \cdot 
					\bigl[  \sed    \sju{1} \cdot [ \sed   \,   \sjd ] \bigr]\bigr] \cdot 
		\bigl[ \sed  \sjd \bigr] \cdot \notag\\
		& \qquad  \bigl[ \sed  \sju{1}  \cdot 
					\bigl[ \sed  \sju{3} \cdot
								\bigl[ \sed \sju{1} \cdot 
									[\sed\,  \sjd
									] 
								\bigr]
					\bigr]  \cdot 
					\bigl[
						\sed    
					\bigr]  
			\bigr]\,.
\end{align}
The  ordered planar tree representation of the recursive $J$ factorisation of the path \eqref{eq:exr1p} (ie.\ \eqref{eq:expp3}) is shown in \figref{fig:rroeeta1}.
%
%======= FIG ===========
%======= FIG ===========
%
The only part of \propref{lem:jfac} that is not obvious is that a $J$ factor is always followed by a $\sjd  $ step (ie.\ the first step to step below the height of a $J$ factor is always a jump down step). 
This can be proved inductively using the level of recursion: At level zero the most general path in $R_1'$ is $\sed \sju{2k+1} \prod_{j=1}^{k+1} \left[ \sed J_{\phi}^\rec{0} \sjd \right]$, thus true. If we assume the proposition is true for all paths containing  level $\ell$ (or smaller) $J$ factors  (illustrated schematically in \figref{fig:Jfac}) then the number of steps in all the  $J$-factors must be  even since, be definition, each starts with a jump step (hence odd) and   each must end on a even down step (ie.\ the step immediately prior to the assumed (odd) $\sjd $ step). 
The number of steps between two consecutive $J$-factors is two (since an odd down must followed by an  even step, hence a down  step). Thus the number of steps in the level $\ell+1$ $J$-factor is even.  
Since it starts with a (odd) jump up step and is even length, it must end with an even (down) step. Thus the next step after the $J$-factor  must be an odd step  and hence a    jump down step. Thus if level $\ell$ is true so is level $\ell+1$ thus, by induction, true for all levels.

\begin{figure}[ht]
	\begin{center}
	\includegraphics[width=30em]{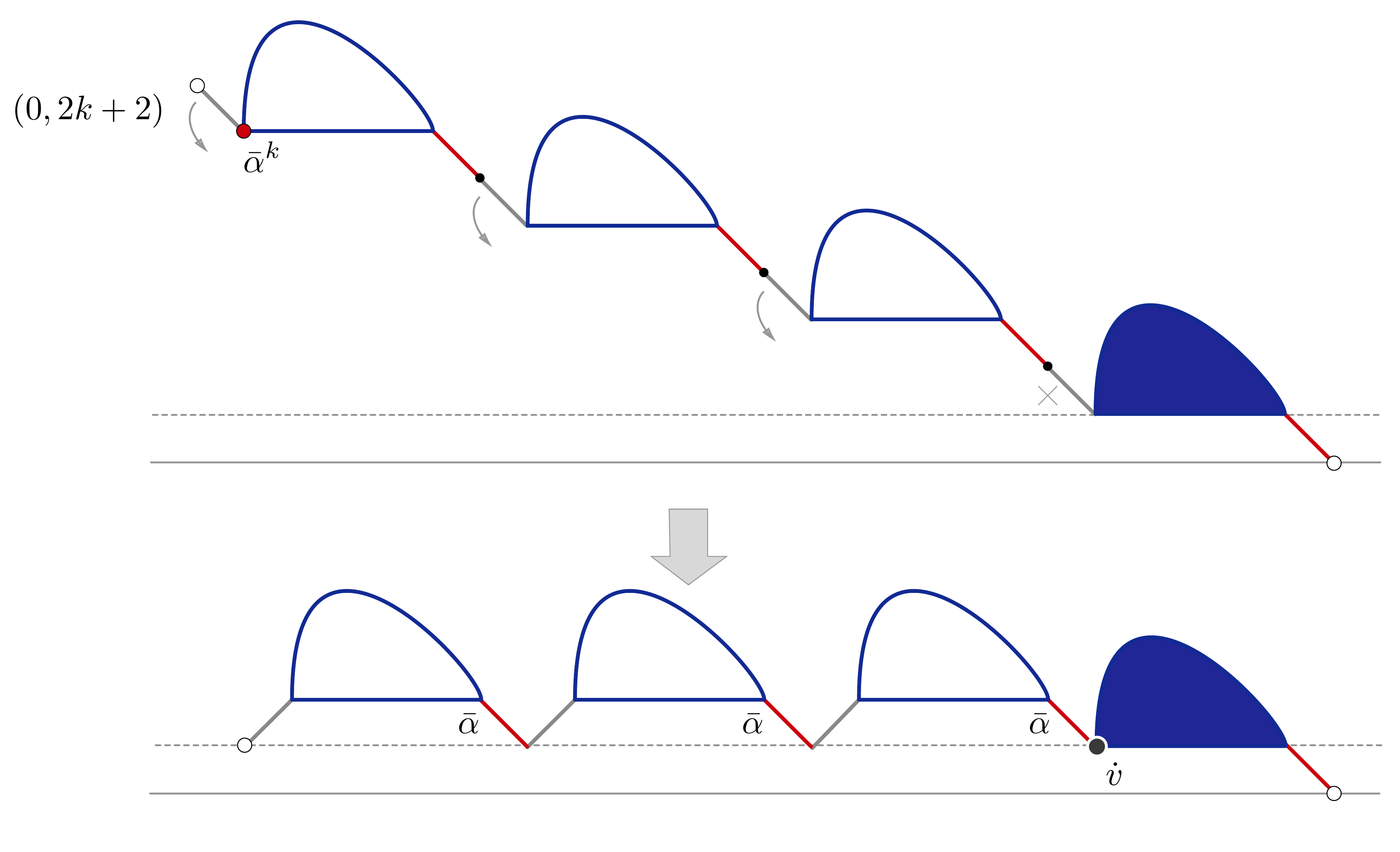}  
	\caption{\it  Schematic representation of the $\Gamma$ map defined by \eqref{eq:gammarotf} whose action gives $R_1^1$ paths according to the $J$ factorisation \eqref{eq:r1paths} -- here $k=3$. The blue $J$ factor contains only $\bb$ weights.}
	\label{fig:rrota1}
	\end{center}
\end{figure}

We now use the $J$-factorisation to biject the paths of $R_1$ to  $R_1^1$ paths which are paths of the form
\begin{equation}\label{eq:r1paths}
  1: \prod_{j=1}^{k} \left[ u \, J_{k_j}\, d \right] \, \dot{v} \,   J_{k_{k+1}}\,  \sjd.  
\end{equation}
All the down, $d$, steps have weight $\ab$ and the up, $u$,  steps have unit weight. The weights of the steps in the $J_{k_{k+1}}$ factor are the same as those of the $R_1$ paths (ie.\ $\bb$ for jump up steps from height one).  The vertex $\dot{v} $ denotes the   vertex which separates the $\ab$ weighted edges from the $\bb$ weighted jump steps in $J_{k_{k+1}}$ and is marked.  The form of the $R_1^1$ paths  are shown schematically in \figref{fig:rrota1} (lower).
%
%======= FIG ===========
%======= FIG ===========
%

\begin{figure}[ht]
	\begin{center}
	\includegraphics[width=30em]{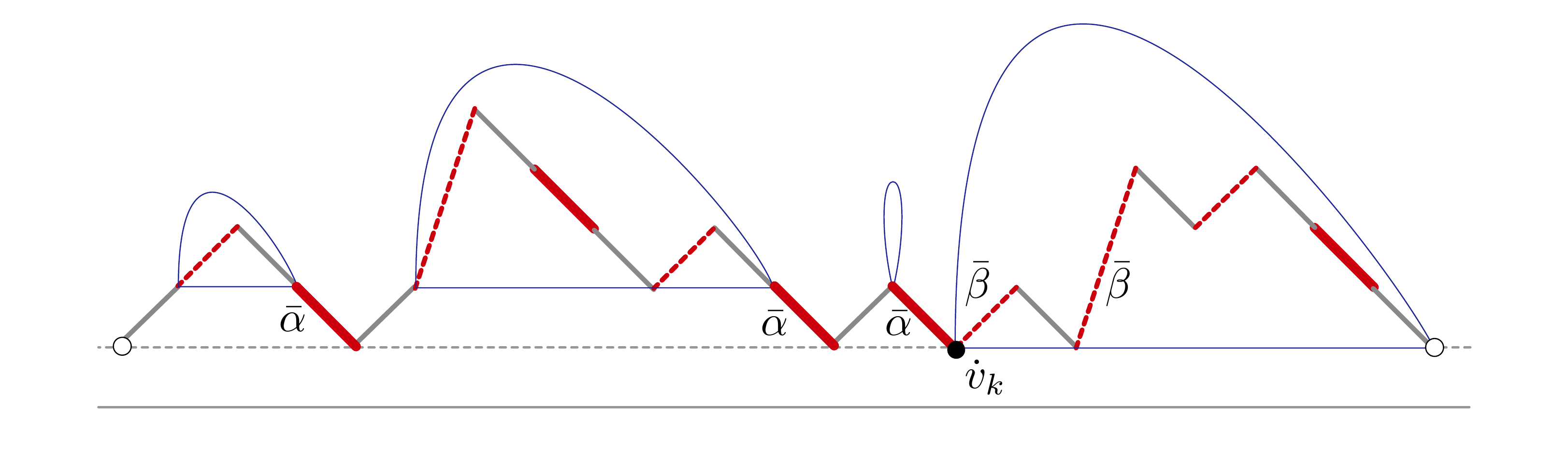}  
	\caption{\it  The result of acting with $\Gamma$ on the example in \figref{fig:r1ex}}
	\label{fig:rrota2}
	\end{center}
\end{figure}

The map $\Gamma: R_1\to R_1^1$, is defined   as follows. If $p \in R_1$ and
\begin{equation}
   \sed\,  p \,\sjd   = 2k+2 :\prod_{j=1}^{k+1} \left[ \sed\, J_{k_j} \,\sjd \right]  
\end{equation}
then the action of $\Gamma$   is defined as
\begin{equation}\label{eq:gammarotf}
	 \Gamma(\sed \, p \,\sjd)=
	\begin{cases}\quad 
		 1:\dot{v} \, p\, \sjd  &\text{if $k=0$}\\
		&\\
		\displaystyle \quad 1: \prod_{j=1}^{k} \left[ u \,J_{k_j} \sjd \,\right] \, \dot{v}  \,
			    J_{k_{k+1}} \sjd   	& \text{if $k>0$}  \\
	\end{cases}
\end{equation}
The weight of each of the explicitly written $\sjd $ steps in \eqref{eq:gammarotf}, to the left of $\dot{v} $ (the marked vertex)  is $\bar{\alpha}$. The height of the  first  and last   vertices of the path $\Gamma(\sed\,  p \,\sjd)$  are the same since  $\Gamma$ has changed $k$ of the $\sed $ down steps of $\sed\,  p\, \sjd$  to up steps and deleted one $\sed $ step -- a height change of the first vertex of of $2k-1$. Since the first vertex of $\sed \, p\, \sjd$ was at height $2k+2$ the net change is to place the vertex at height one.  
This map is illustrated schematically in \figref{fig:rrota1} and for a particular example in \figref{fig:rrota2}. It is straightforward to show  $\Gamma$ is a bijection and so we omit the details.
%
%======= FIG ===========
%======= FIG ===========

\noindent\emph{Stage 2:} The $R_1^1$ paths have the simple factored form given by \eqref{eq:gammarotf} which we  now biject     to $R_4$  paths by acting independently on each of the $J$-factors in \eqref{eq:gammarotf} to produce $D$-factors. The action of the map $\Gamma'$ on $p\in R_1^1$   is given in terms of the form \eqref{eq:r1paths} as
\begin{equation}\label{eq:gponr11}
 \Gamma'(   p  ) =  
       \prod_{j=1}^{k} \Gamma'\left( u \, J_{k_j}\, d \right) \, \dot{v} \,  \Gamma'\left(  J_{k_{k+1}}\right)\,  \sjd.  
\end{equation}
and the action of $\Gamma'$ on a $u \, J_{k}\, d $ factor  is defined   recursively using the factorisation \eqref{eq:jfac2} (omitting the level superscripts) by
\begin{equation}\label{eq:r14biject}
	 \Gamma' \left( \sed J_k\, \sjd \right) =
	\begin{cases}
		\displaystyle	\sdd 								& \text{if $J_k=\phi$}\\
		\displaystyle	\sdd \sdu \cdot \sdu^{k } \prod_{j=1}^{k+1} \Gamma '\left(  \sed J_{k_j} \sjd \right)  
						 	& \text{if $J_k\ne\phi$}
	\end{cases}
\end{equation}
Thus $\Gamma'$ has replaced the first $\sed$  of \eqref{eq:jfac2} (of the righthand side case two)  by $d$ and the $u_{2k+1}$ step by $u^{k+1}$. Any $\ab$ weighted $d$ step retains the $\ab$ weight under the action of $\Gamma'$. All the $\bb$ weights are  associated with the jump up steps (from height one -- see \figref{fig:r1ex}) in the rightmost $J$-factor ie.\ $J_{k_{k+1}}$ and under $\Gamma'$ the $\bb$ weight is associated with the leftmost $u$ step of \eqref{eq:r14biject}. 

We define $D$ by
\begin{equation}\label{eq_gpdfact}
\sdd D  = \Gamma' \left( \sed J_k\, \sjd \right).
\end{equation}
where the use of $D$ signifies that $\Gamma'$ produces elevated Dyck subpaths (proved below).

\begin{figure}[ht]
	\begin{center}
	\includegraphics[width=30em]{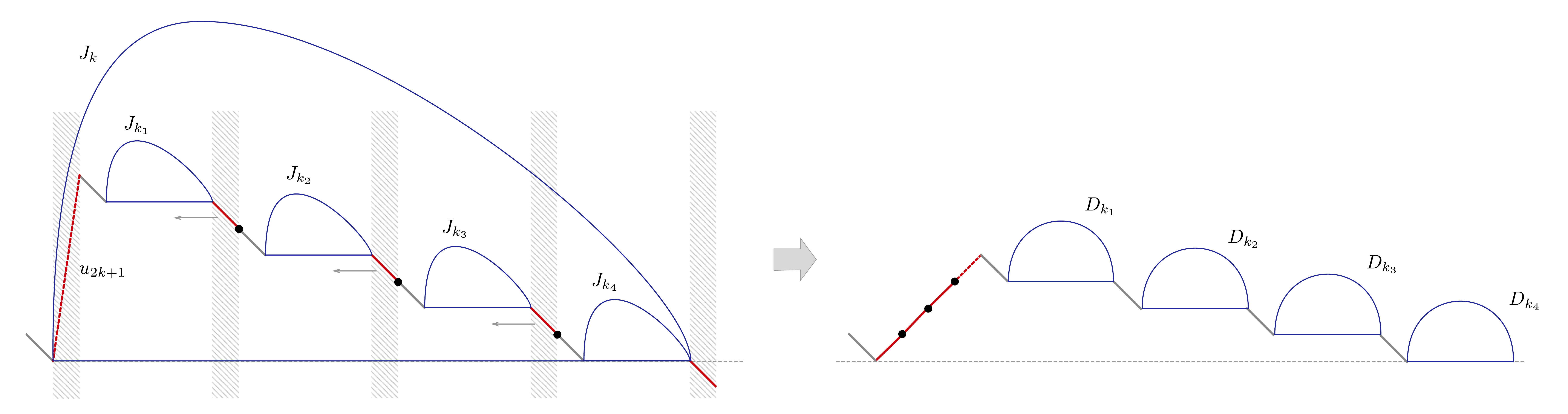}  
	\caption{\it  Schematic representation of the $\Gamma'$ map defined by \eqref{eq:r14biject} -- (here $k=3$) giving a Dyck path.}
	\label{fig:jtodmap}
	\end{center}
\end{figure}

\begin{figure}[ht]
	\begin{center}
	\includegraphics[width=30em]{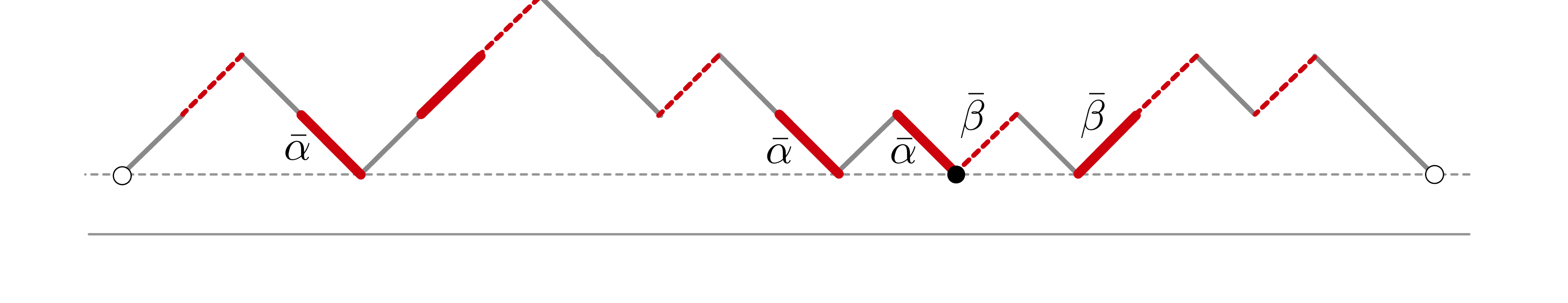}  
	\caption{\it  The result of acting with $\Gamma'$ on the example in \figref{fig:rrota1}}
	\label{fig:r1ex1_r4}
	\end{center}
\end{figure}

The $\Gamma'$ map is illustrated schematically in \figref{fig:jtodmap}.
%======= FIG ===========
%======= FIG ===========
For example, with $\Gamma'$ applied to \eqref{eq:exr1p} via the factorisation \eqref{eq:jfacex} the image path is:
\begin{align}
		 \Gamma'(\sed J_3\, \sjd) & =  
		  \sdd  \sdu^4\cdot  \Gamma'( \sed J_0\, \sjd ) \cdot \Gamma'(\sed J_1\, \sjd ) \cdot \Gamma'(\sed  \sjd )  \cdot \Gamma'( \sed J_0\, \sjd ) \notag \\	
		  & =  
		  \sdd  \sdu^4\cdot  [ \sdd  \sdu^1 \cdot \Gamma'( \sed  \sjd )] \cdot [ \sdd  \sdu^2 \cdot  \Gamma'( \sed  \sjd   ) \cdot 
					\Gamma'(  \sed   J_0  \,   \sjd  )] 
		   \cdot \sdd   \cdot 
		   \dot{v}_3 \cdot 
		   [\sdd  \sdu^1 \cdot \Gamma'(   \sed    J_1\, \sjd ) \bigr])] \notag \\	
		&= \sdd  \sdu^4\cdot  [ \sdd  \sdu^1 \cdot  \sdd   ] \cdot [ \sdd  \sdu^2 \cdot   \sdd  ] \cdot 
			[\sdd  \sdu^1 \cdot \Gamma'(  \sed     \sjd  )] 
		   \cdot \sdd   \cdot \notag\\
		   &\qquad 
		   [\sdd  \sdu^1 \cdot 
		    \sdd  \sdu^2 \cdot  \Gamma'( \sed  J_0\, \sjd   )\cdot \Gamma'(  \sed      \sjd   ) ] \notag \\
		 &= \sdd  \sdu^4\cdot  [ \sdd  \sdu^1 \cdot  \sdd   ] \cdot [ \sdd  \sdu^2 \cdot   \sdd  ] \cdot 
			[\sdd  \sdu^1 \cdot   \sdd      ] 
		   \cdot \sdd    \cdot \notag\\
		   &\qquad 
		   [\sdd  \sdu^1 \cdot 
		    \sdd  \sdu^2 \cdot  \sdd  \sdu^1 \cdot  \Gamma'( \sed  \sjd   )\cdot  [ \sdd     ] ] \notag \\
		 &= \sdd  \sdu^4\cdot  [ \sdd  \sdu^1 \cdot  \sdd   ] \cdot [ \sdd  \sdu^2 \cdot   \sdd  ] \cdot 
			[\sdd  \sdu^1 \cdot   \sdd      ] 
		   \cdot \sdd    \cdot \notag\\
		   &\qquad 
		   [\sdd  \sdu^1 \cdot 
		    \sdd  \sdu^2 \cdot  \sdd  \sdu^1 \cdot  [ \sdd  ]\cdot  [ \sdd     ] ] \notag \\
		& =  \sdd  \sdu^4 \sdd  \sdu^1  \sdd   \sdd  \sdu^2   \sdd   \sdd  \sdu^1    \sdd    \sdd   \cdot   \sdd  \sdu^1  
		    \sdd  \sdu^2   \sdd  \sdu^1   \sdd   \sdd    
\label{eq:jfacex}
\end{align}
The result of acting with $\Gamma'$ on the example in \figref{fig:rrota2} is shown in figure \figref{fig:r1ex1_r4}.
%======= FIG ===========
%======= FIG ===========
The path configurations after acting with $\Gamma'$ is that of $R_4$ except for that the $\bb$ weights are on the up step (from height one to two) rather than on the down step (form height two to one) however this is readily fixed just by moving the weight across.

We now prove by induction on the level of recursion  that the $D$ factor of  \eqref{eq_gpdfact} is an elevated Dyck path.
Clearly the step set of $D$ is that of Dyck paths. What needs justification is that that paths in $D$ start and end at the same height and no vertices of the path are below that of the initial vertex. 

Re-instating  the level of recursion with a superscript and subscripts to distinguish the $D$ factors, the initial step of the induction corresponds to  with   $J^{(0)}_k=\phi$ in which case $d\, D^{(0)}_k=\Gamma \left( \sed  \sjd \right)=d$, thus $D^{(0)}_k=\phi$ which is an (empty) Dyck path. 
Inducting from level $\ell$ to $\ell+1$ we have
\begin{align}
d\, D^{(\ell+1)}_k  =\Gamma' \left( \sed J^{(\ell+1)}_k\, \sjd \right) 
	&=\sdd \sdu^{k+1} \prod_{j=1}^{k+1} \Gamma '\left(  \sed J^{(i_j)}_{k_j} \sjd \right)\\
	&=\sdd \sdu^{k+1} \prod_{j=1}^{k+1} \left[ d\, D^{(i_j)}_{k_j}   \right]
\intertext{where, as in \eqref{eq:eq_lamx},  $\ell=\max\set{i_j \st j=1...k+1}$, thus}
D^{(\ell+1)}_k 
	&=   \sdu^{k+1} \prod_{j=1}^{k+1} \left[d\,  D^{(i_j)}_{k_j}   \right].   \label{eq_repfg3}
\end{align}
If we assume for all levels $i_j\le \ell$  each $D^{(i_j)}_\ell$    is an elevated Dyck path (and hence the first and last vertices are the same height) and since the prefix $\sdu^{k+1} $ in \eqref{eq_repfg3}, goes up $k$ steps and the product   steps down $k+1$ times (ie.\ the $k+1$, $d$ steps), the righthand side is also a Dyck path, that  is $D^{(\ell+1)}_k$ is a Dyck path, thus by induction the proposition is true.

%! R_2 -> R_4 =========
\subsection{Proof of Equivalence of the $R_2$ and $R_4$ path representations}
%!***************************************************************

%

We prove this equivalence in four stages.  The four stages are connected by either a bijection or a sign reversing involution. The five intermediate  sets of paths involved, $R_2^i$, $i=1..5$ are defined when each stage is discussed in detail  below.
\begin{itemize}

\item[Stage 1.] $R_2\,\xrightarrow{\,\ka^2\,}\, R_2^1 \, \xrightarrow{\,\Phi_2^{12}\,}\,  R_2^2$. 
A sign reversing  involution, $\Phi_2^{12}$, which reduces the infinite sum \eqref{eq:lpf2} over $R_2$ paths to a 
finite sum   over $R_2^2$ paths. The involution acts on an enlarged path set $R_2^1$, obtained from  $R_2$ paths by expanding $\ka^2=1-cd$. The fixed point set of $\Phi^{12}_2$ is the set of $R_2^2$ paths.

\item[Stage 2.] $R_2^2 \,\xrightarrow{\,\Gamma_2^{23}\,} \, R_2^3$. 
%The sum over $R_2^1$ paths is replaced by a sum over $R_2^{2}$  paths  using a bijection, 
The bijection  $\Gamma_2^{23}$  `pulls down' the first and last vertices of each path thus  replacing the sum over $R_2^2$ paths  by a sum over $R_2^{3}$ paths (which start and end at height one).

\item[Stage 3.] $R_2^3\,\xrightarrow{\,\Gamma_2^{34}\,}\, R_2^4$.  
The bijection  $\Gamma_2^{34}$  `lifts' the $R_2^3$ paths above the surface to give $R_2^4$ paths (which have no height zero vertices).

\item[Stage 4.] $R_2^4 \,\xrightarrow{\,c,d\to \ab,\bb}\, R_2^5 \,\xrightarrow{\,\Phi_2^{56}\,}\,  R_4$. 
The final sign reversing involution, $\Phi_2^{56}$, replaces the $c$ and $d$ weighted paths of $R_2^4$ with   $\ab$ and $\bb$ weighted paths. The involution acts on an enlarged set of paths, $R_2^5$,  obtained by expanding $c=1-\ab$ and $d=1-\bb$.   The fixed point set is the path  set  $R_4$.

\end{itemize}
In summary,
\begin{equation}\label{eq:r2r4map}
 R_2\,\xrightarrow{\,\ka^2\to1-cd\,}\, R_2^1 \,
\xrightarrow{\,\Phi_2^{12}\,}\,  R_2^2 \,
\xrightarrow{\,\Gamma_2^{23}\,} \, R_2^3\,
\xrightarrow{\,\Gamma_2^{34}\,}\, R_2^4 \,
\xrightarrow{\,c,d\to \ab,\bb}\, R_2^5 \,
\xrightarrow{\,\Phi_2^{56}\,}\,  R_4
\end{equation}
 We now expand on each of the four stages.

%! %%%%%%%% Stage 1 %%%%%%%%%
\paragraph{Stage 1.} $ R_2\,\xrightarrow{\,\ka^2\,}\, R_2^1 \,\xrightarrow{\,\Phi_2^{12}\,}\,  R_2^2 $.
The sign reversing involution is defined on the set of paths  $ R_2^1 $ which is constructed   by using $\ka^{2}=1-cd$ to enlarge the size of the weighted set $R_2$ (which has weights given by \eqref{eqs:w2}). 
Thus for each weighted path $\om\in R_2$ (which always has a factor of $\ka^{2}$ in its weight) we replace by two paths $\om_{1}$ and $\om_{2}$, where $\om_{1}$ is the same sequence of steps as $\om$, but the initial and final vertex weights are $w^{i}((0,2k+1))=c^{k} $ and $w^{f}((2L,2k'+1))=d^{k'} $ (ie.\ no factors of $\ka$). 
Similarly, $\om_{2}$ is the same sequence of steps as $\om$, but the initial and final vertex weights are $w^{i}((0,2k+1))=-c^{k+1} $ and $w^{f}((2L,2k'+1))=d^{k'+1}$ ie.\ each vertex has an extra factor of $c$ (or $d$), and an overall negative weight).  Thus we have that
\begin{align}
		Z_{2L}^{(2)}&=\sum_{\om\in 
R_2^1}W^{(2)}_2(\om) 
\end{align}
where the weight $W^{(2)}_2 $ is as just  explained.   The   $R_2^2$ paths are a subset of the $R_2$ paths, given by
\begin{equation}\label{eq:r2ppaths}
 R_2^2=\set{p\in R_2 \st \text{$p$ has at least one vertex of height one}}
\end{equation}

We will now show that $R_2^2$ is the fixed point set of $R_2^1$ under the sign reversing involution  $\Phi_2^{12}$  defined below. The signed set 
  $\Om^{(2)}=R_2^1=\Om^{(2)}_{+}\cup \Om^{(2)}_{-}$   is defined by: 
\begin{align}
		\Om^{(2)}_{+}&=\{\om\suchthat 
		\text{$\om\in R_2^2$ and $W^{(1)}_2 (\om)>0$} \}\\
		\Om^{(2)}_{-}&=\{\om\suchthat 
		\text{$\om\in R_2^2$ and $W^{(1)}_2 (\om)<0$}\}.
\end{align}
The involution $\Phi_2^{12}\,:\,\Om^{(2)} \to \Om^{(2)}$ is defined by three cases.  Let
$\om \in \Om^{(2)} $, $\om'=\Phi_2^{12}(\om)$ and let $v_0$ be the first vertex of $\om$ and $v_{2L}$ the last.  Recall, $w(v)$ is the weight of vertex $v$.
\begin{itemize}

\item[Case 1.]\textbf{(Negative weight.)}
  If   $v_{0}=(0,2k+1)$, $v_{2L}=(2L,2k'+1)$, $k,k'\ge 0$ and $w^{i}(v_{0})=-c^{k+1} $ then $\om'$ is a path with the same sequence of steps as $\om$, but initial vertex $v_{0}'=(0,2k+3)$, final vertex $v_{2L}=(2L,2k'+3)$ (ie.\ is $\om$ `pushed up' two units), and has vertex weights $w (v'_{0})=c^{k+1} $ and $w (v'_{2L})=d^{k+1}$.
For any $\om$, $\om'$ always exists  and has opposite sign to $\om$, thus $\Phi_2^{12}$ is sign reversing for this case.

\item[Case 2.]  \textbf{(Positive weight, no height one vertices.) }
If    $v_{0}=(0,2k+1)$, $v_{2L}=(2L,2k'+1)$, $k,k'\ge 1$, $w (v_{0})=c^{k} $ and $\om$ has \textit{no vertex with height one}, then $\om'$ is a path with the same sequence of steps as $\om$, but initial vertex $v_{0}'=(0,2k-1)$, final vertex $v_{2L}=(2L,2k'-1)$ (ie.\ is $\om$ ``pushed down'' two units), and has vertex weights $w (v'_{0})=c^{k} $ and $w (v'_{2L})=d^{k} $.
Since $\om$ no height one vertices, all its vertices have height greater than two, thus when $\om$ is pushed down no vertices have height less than zero and hence $\om'\in \Om^{(2)}$.  For any $\om$ in this case, $\om'$ always exists  and has opposite sign to $\om$, thus $\Phi_2^{12}$ is sign reversing for this case.

\item[Case 3.]  \textbf{(Positive weight, at least one height one vertex.)} If $\om$ has positive weight and at least one vertex with height one, then $\om'=\om$.
		 
\end{itemize}

Clearly, if $\om$ corresponds to Case 1, then $\om'$ is a unique path
corresponding to Case 2 and visa versa.  Case 3 is the fixed point set. Since the fixed point set   paths are in the positive set, $\Om_+^{(2)}$, they have weight $c^k$ for the initial, height $2k+1$ vertex and weight $d^{k'}$ for the last, height $2k'+1$ vertex.
Thus $\Phi_2^{12}$ is a sign reversing involution with fixed point set the subset of $R_2^1$ paths with at least one vertex at height one and positive weight  ie.\ $R_2^2$ paths. 
	
The paths in $R_2^3$ have at least one vertex with height one and may have many with height zero. We `biject away' the latter subset in the next stage.
%The paths of $R_2^3$ can start and end at different (odd) heights. We now use a bijection to a set of paths that start and end at height one.

%! %%%%%%%%% Stage 2 %%%%%%%%%

\paragraph{Stage 2.} $R_2^2 \,\xrightarrow{\,\Gamma_2^{23}\,} \, R_2^3\,$.
We now map the path set to a subset $R_2^3$ of $R_2^2$ paths which do not intersect the line $y=0$. In order to do this the resulting paths have to carry a ``dividing'' line (or equivalently a marked vertex). Thus, if
\begin{align} \label{eq:r23paths}
 \hat{R}_2^2&=\set{p \in R_2^2 \st \text{$p $ has no height  zero vertices}}
\intertext{then}
 {R}_2^3&=\ \text{set of paths of $\hat{R}_2^2$  with one   height one vertex marked.} 
\end{align}
That is, if  $p \in {R}_2^3$ has  $m$ vertices with height one, then $p $ produces $m$ paths in ${R}_2^3$ each one with one of the $m$  vertices marked.

Let $p \in R_2^2$. If $p $ starts at height $2k+1$ and ends at height $2\ell+1$ then, using a similar factorisation to the $D$-factorisation of the $R_1$ to $R_4$ bijection -- \lemref{lem:catfac},   $p $ can be factorised as 
\begin{equation}\label{eq:g22fact}
 \left[\prod_{n=1}^{2k} D_n d  \right] B \left[ \prod_{m=1}^{2\ell}   u D'_m\right]
\end{equation}
where $D_n$ and $D'_m$ are (possibly empty ) elevated Dyck paths, $u$ an up step, $d$ a down step and $B$ is defined by the fact that $uBd$ is  a Dyck path. 
%
%>>>>>>>>>>>>>>>>>>>>>>>>>>>>>>>>>>>>>>>>>>>>>>>>>
\begin{figure}[ht!]
	\centering
\includegraphics[width=30em]{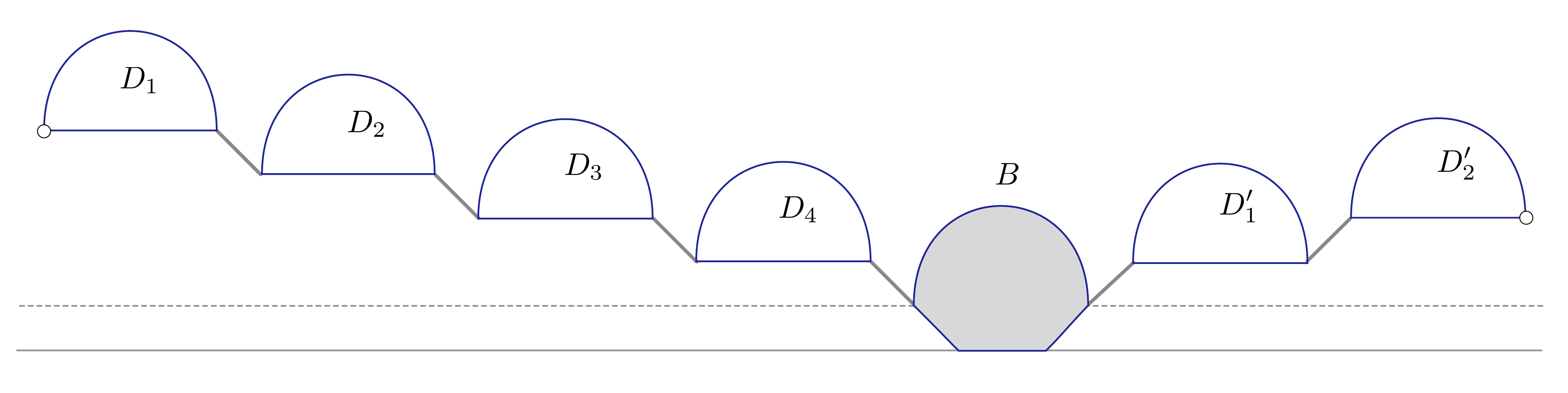}
	\caption{\it An schematic representation of the factorisation \eqref{eq:g22fact} (here $k=2$, $\ell=1$). }
\label{fig:gb2}
\end{figure}
%>>>>>>>>>>>>>>>>>>>>>>>>>>>>>>>>>>>>>>>>>>>>>>>>>
%
That is, $B$ is the subpath of $p $ which  is made of only up and down steps and whose first vertex is the leftmost height one vertex of $p $ and whose last vertex is the rightmost height one vertex of $p $. If $k$ or $\ell$ is zero then the respective product  is absent. The factorisation is shown schematically in \figref{fig:gb2}.

\begin{figure}[ht!]
	\centering
\includegraphics[width=30em]{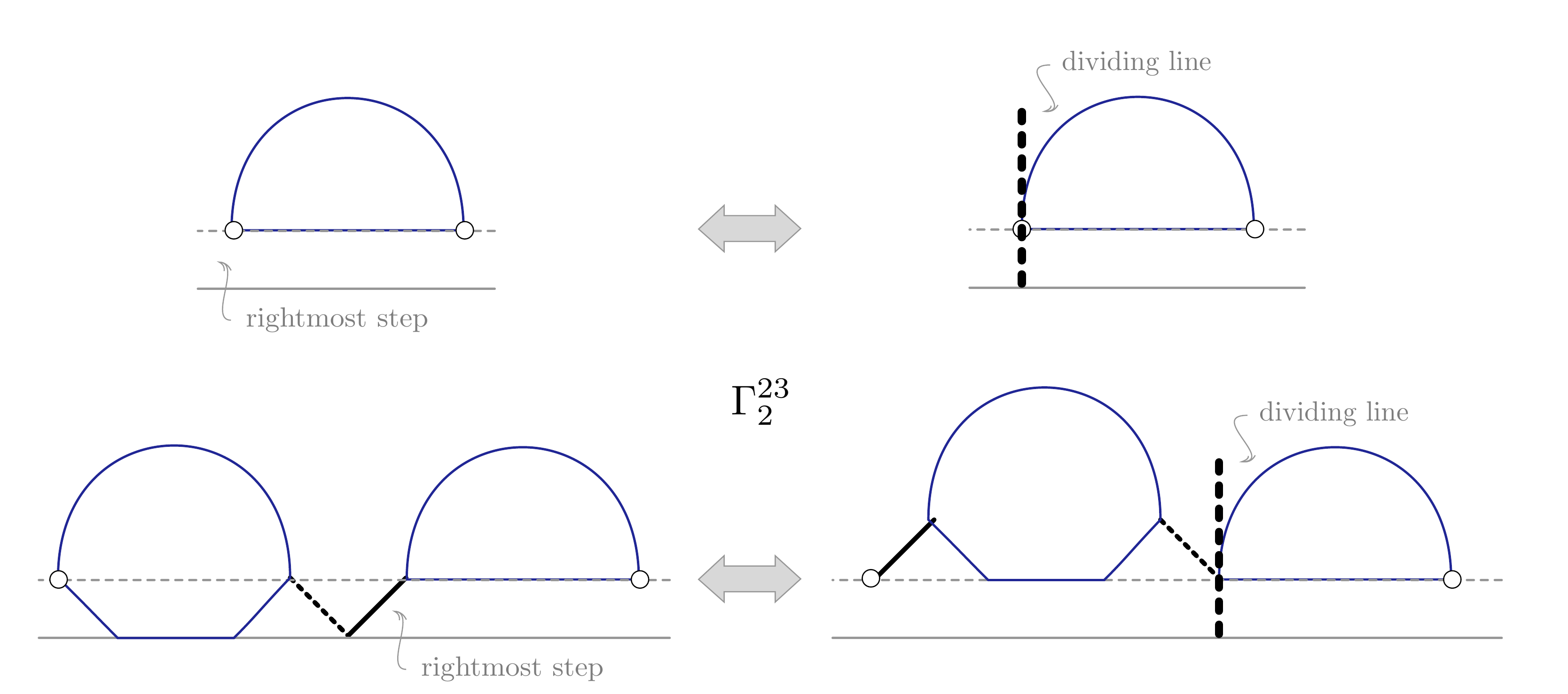}
	\caption{\it An schematic representation of the action of the  bijection 
$\Gamma_2^{23}$ on the $B$ factor at stage 2.  Case with no  step below $y=1$ (upper) and the case with at least one step below $y=1$ (lower). }
\label{fig:gb}
\end{figure}

We now construct a map, $\Gamma_2^{23}: R_2^2\to R_2^3$, that eliminates all steps  of the subpath $B$, below $y=1$ and replaces it with a path, $\hat{B_L}| \hat{B_R}$, which has no height zero vertices but has a `dividing line' (or marked vertex)
-- see \figref{fig:gb}.     
%
%>>>>>>>>>>>>>>>>>>>>>>>>>>>>>>>>>>>>>>>>>>>>>>>>>
%>>>>>>>>>>>>>>>>>>>>>>>>>>>>>>>>>>>>>>>>>>>>>>>>>
%

The map  $\Gamma_2^{23}$ acts on $B$ as follows. If $B$ has no steps below $y=1$, then $ B'=\Gamma_2^{23}(B)=| B$, where $ | $ denotes a vertical dividing line drawn through the leftmost vertex of $B$. If $B$ has at least one step below $y=1$, then let $u'$ be the rightmost (up) step from $y=0$ to $y=1$. 
Thus $B$ factorises as $B=w_1 u' w_2$, and then $B'=u' w_1 |  w_2 $, where $ | $ denotes a vertical dividing line drawn through the vertex between $w_1$ and $w_2 $. 
Note, since $u'$ is an up step, none of the steps of the subpath $u' w_1$ intersect $y=0$. Thus $B'$ does not intersect $y=0$. 
%Hence $B'=B_L | B_R\in R_2^3$. 
The map $\Gamma_2^{23}$ acting on all factors of the form of $B$ is readily seen to be injective and surjective and thus a bijection (the dividing line shows where the first up step has to be moved under the action of the inverse map $\left(\Gamma_2^{23}\right)^{-1}$).

The action of $\Gamma_2^{23}$ on $p \in R_2^2$ only depends on its  $B$ factor and is defined as
\begin{equation}\label{eq:g212act}
 \Gamma_2^{23}(p )=\Gamma_2^{23}\left( \prod_{n=1}^{2k} D_n d   \cdot B \cdot \prod_{m=1}^{2\ell}   u D'_m \right)
= \prod_{n=1}^{2k} D_n d \cdot\Gamma_2^{23}(B) \cdot \prod_{m=1}^{2\ell}   u D'_m 
\end{equation}
with the weight of all vertices unchanged. Thus the path $\Gamma_2^{23}(p )$ has the  same weight as $p $, does not intersect $y=0$ and has a dividing line, that is, $\Gamma_2^{23}(p )\in R_2^3$.

%! %%%%%%%%% Stage 3 %%%%%%%%% 

\begin{figure}[ht!]
	\centering
\includegraphics[width=30em]{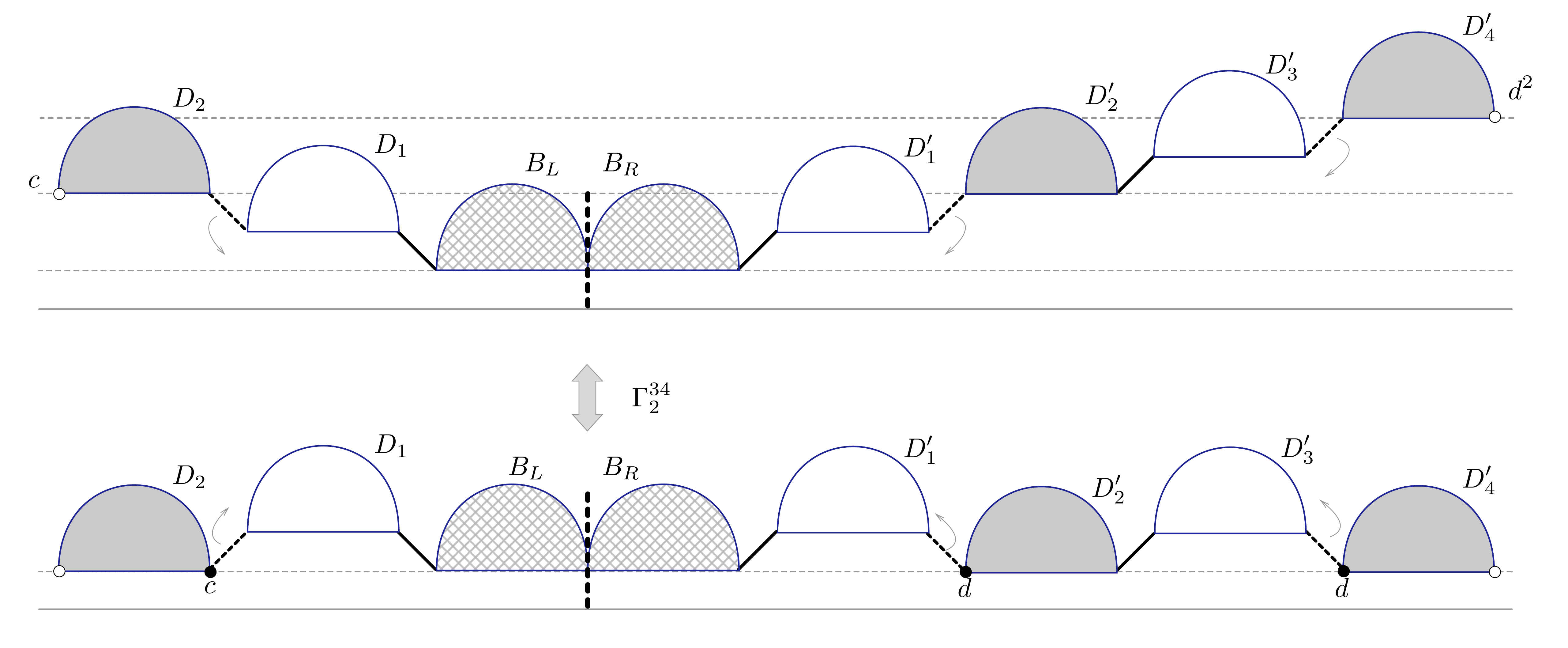} 
	\caption{\it An schematic representation of the bijection 
$\Gamma_2^{34}$ of stage 3.}
\label{fig:gr}
\end{figure}

 \paragraph{Stage 3.} $R_2^3\,\xrightarrow{\,\Gamma_2^{34}\,}\, R_2^4 $. 
The map  $\Gamma_2^{34} $  `rotates down'  the initial and final vertices of the path to produce a path which starts and ends at $y=1$, but has a subset of ``marked'' $c$ and $d$ height one vertices. This is a simple extension of the same map given in \cite{brak:2001yf} and hence we only discuss it briefly here.  It is illustrated schematically in \figref{fig:gr}. 
%>>>>>>>>>>>>>>>>>>>>>>>>>>>>>>>>>>>>>>>>>>>>>>>>>
%>>>>>>>>>>>>>>>>>>>>>>>>>>>>>>>>>>>>>>>>>>>>>>>>>

Let $p \in R_2^3 $ start at $y=2k+1$,  and end at $2\ell+1$,  (and 
hence has weight $c^k d^\ell$). Using the factorisation \eqref{eq:g212act},
\begin{equation}\label{eq:eq:g23fact}
	p = \prod_{n=1}^{2k} D_{2k-n+1} d \cdot B_L | B_R \cdot \prod_{m=1}^{2\ell}   u D'_m 
\end{equation}
we can define $\Gamma_2^{34}$ by 
\begin{equation}\label{eq:g23ase3}
	\Gamma_2^{34}(p )=\left( 
	\prod_{n=1}^{k} D_{2k-2n+2} \cdot\dot{v}(c) \cdot u \, D_{2k-2n+1} d\right) 
	\cdot B_L | B_R \cdot 
	\left( \prod_{m=1}^{\ell }   u \,D'_{2m-1} d \cdot\dot{v}(d) \cdot   D'_{2m} \right)
\end{equation}
  where, $\dot{v}(c)$ and $\dot{v}(d)$ represent  a 
 marked  vertex between the two steps where it occurs (and is weighted $c$ and $d$ respectively) -- see \figref{fig:gr}. Each mark  to the left of the dividing line carries weight $c$ and each of those to the right of $|$ carry a weight $d$.  The inverse map $\left(\Gamma_2^{34}\right)^{-1}$ uses the marked vertices to fix the step change $u\to d$.

%%%%%%%%%%: Stage 4 %%%%%%%%%

\paragraph{Stage 4.} $R_2^4 \,\xrightarrow{\,c,d\to \ab,\bb}\, R_2^5 \,\xrightarrow{\,\Phi_2^{56}\,}\,  R_4$. 
In the final stage we define an involution, $\Phi_2^{56}$ whose fixed point set is $R_4$  with weights $\ab$ and $\bb$. Starting with the set of all paths given at Stage 3, ie.\ of the form\eqref{eq:g23ase3}, we construct a larger set of paths, $R_2^5$, using the same construction of Stage 1, that is,  by replacing all weights  $c$ with $\ab-1$ and all weights $d$ by $\bb-1$. 
Thus each path in $R_2^4$,  which has weight $c^k d^\ell$, maps to $2^{k+\ell}$ paths. Combinatorially, this set has all marked vertices, $\dot{v}(c)$ (ie.\ to the left of $|$), replaced by either a weight of $-1$ or $\ab$ and all marked vertices, $\dot{v}(d)$ (ie.\ to the right of $|$), replaced by either a weight of $-1$ or $\bb$. 
All remaining vertices of the path intersecting $y=1$, \textit{except} that intersecting the dividing line, will be labeled with `$+1$'.  The weight of a given path is a product of all the $\ab$, $\bb$ and $-1$ factors. Thus the weight of the path will be negative if there are an odd number of factors of $-1$. 

This construction defines the elements of the set $\Omega=\Om_+\cup\Om_-$ where $\Om_+$ contains the positive weighted paths and $\Om_-$ the negative weighted paths.
 \begin{figure}[ht!]
 	\centering
\includegraphics[width=30em]{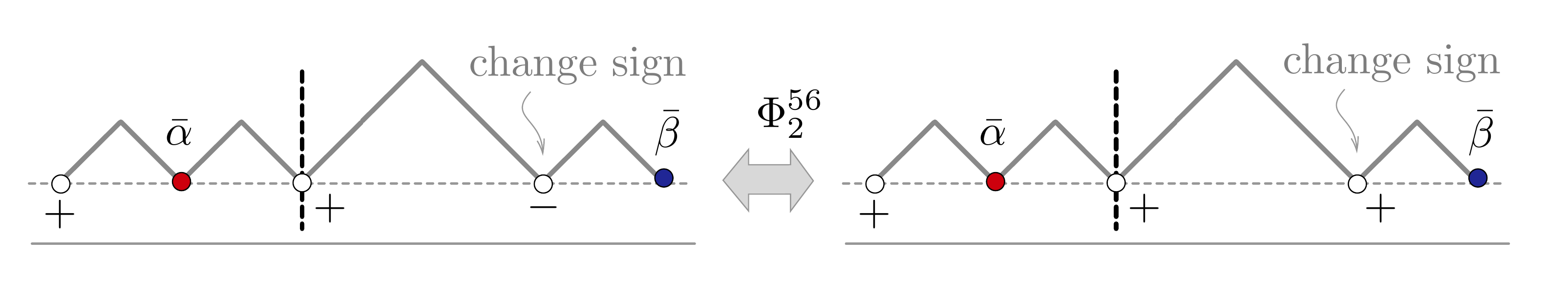} 
 	\caption{\it  An example of a cancelling a pair of partially marked paths -- the rightmost $\pm 1$ weighted vertex is changed to a $\mp 1$ weight. }
 \label{fig:gc4s}
 \end{figure}
The involution, $\Phi_2^{56}$, is straightforward:
If $p \in\Omega$ has no $+1$ or $-1$ vertices then $p '=\Phi_2^{56}(p )=p $. All these cases obviously have positive weight with all height one vertices  to the left of the dividing line carrying weight $\ab$ and those to the right, weight $\bb$. These are the fixed point paths and are clearly $R_4$ paths (after deleting the dividing line -- which is no longer necessary). 
If $p \in\Omega$ has at least one $+1$ or $-1$ vertex then $p '=\Phi_2^{56}(p )$ is the same weighted path as $p $ except the rightmost signed vertex has opposite sign (and hence $p '$ has the same weight as $p $ except of opposite sign) -- see \figref{fig:gc4s}.  Clearly, $\left(\Phi_2^{56}\right)^2=1$.

%  \begin{figure}[ht]
% \begin{center}
% %\includegraphics{\figPathSketch{UoutlineTest.pdf}
% \includegraphics{intercolate.pdf}
% \caption{ $2$-elevated Dyck path (left) and bijected jump step subpath (right) }
% \label{fig:js}
% \end{center}
% \end{figure}

%An example for each of the four stages is shown in figure 
%\ref{fig:exewr}.

%! R_3 -> R_4 ============
\subsection{Proof of Equivalence of the $R_3$ and $R_4$ path representations}
%***************************************************************

\begin{figure}[ht!]
	\centering
	\includegraphics{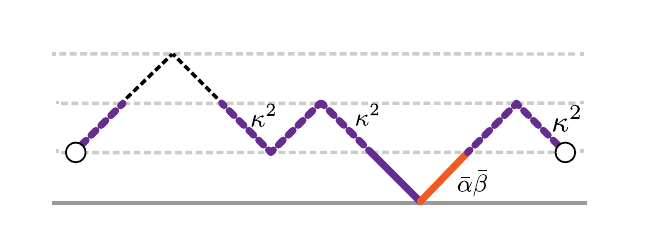}
	\caption{\it  An example of a re-weighted $R_3$ path.}
	\label{fig:r3ppaths}
\end{figure}

%
We prove the equivalence using a sign reversing involution, $\Phi_3 $.  The fixed point set will be the set of paths $R_4$.  Before defining the signed set of the involution we re-weight  the steps of the $R_3$ paths as follows. The paths in $R_3$ have steps from height two to one and height one to two each weighted by $\kappa$ (see \defref{def:r2p}). Since all the paths in $R_3$ start and end at height one, all paths have an even number of steps between heights two and one and thus each path has an even degree $\kappa$ weight ie.\ $\kappa^{2k}$ (readily proved by induction on the length of the path). Thus rather than have $\kappa$ weights associated with up and down steps we associate a $\kappa^2$ weight only with a down step (from height two to one).  Similarly   there are an even number of  steps between heights zero and one. These carry weights $\ab$ and $\bb$ so we collect the two weights together to form a single $\ab\bb$ weight associated with the up step from height zero to one and give the down step unit weight. Call this reweighed path set, $R_3'$. An example is shown in \figref{fig:r3ppaths} (which is a re-weighting of the example in \figref{fig:state}).

We now  increase the size of $R_3'$ by  expanding all $\ka^{2}=\ab+\bb-\ab\bb$ weights.  
Thus any path, $\om$ with an edge, $e_{n}$ with weight $\ka^{2}$ gives rise to three paths, $\om_{1}$, $\om_{2}$ and $\om_{3}$, with the same step sequence, but different weights:  $\om_{1}$ is the same path as $\om$, but edge $e_{n}$ has weight $\ab$.  Similarly, for $\om_{2}$, edge $e_{n}$ has weight $\bb$ and for $\om_{3}$, edge $e_{n}$ has negative weight $-\ab\bb$. Thus if the path has a weight factor $\kappa^{2k}$ it will give rise to $3^{k}$ paths. Call this expanded set, $R_3^2$. Note, all the $-\ab\bb$ weights are between heights two and one whilst all the $\ab\bb$ weights are between heights zero and one.
\begin{figure}[ht!]
	\centering
	\includegraphics[width=30em]{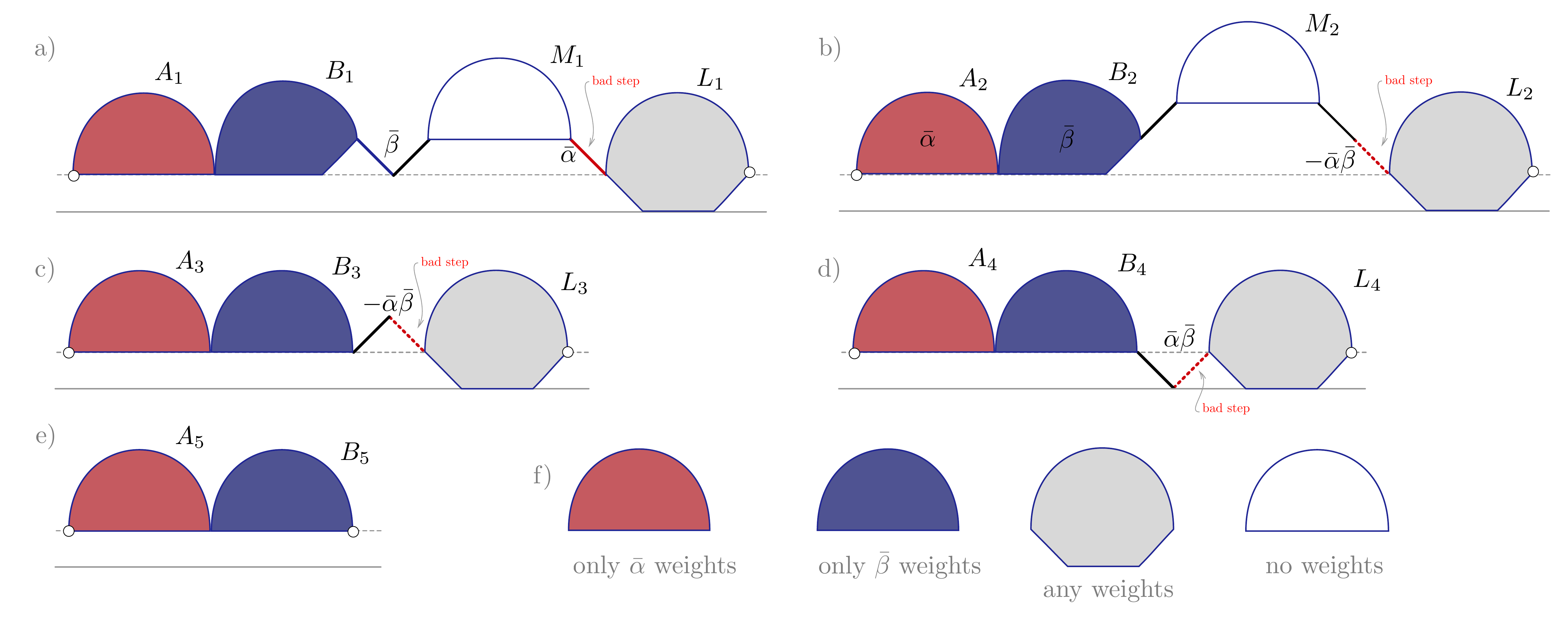}
	\caption{\it a) -- e) Schematic representations of the five possible factorisations of $R_3$ paths  as defined by the position of the `bad' step -- see \eqref{eq:f5}. f) Key for weight structure of subpath factors.}
	\label{fig:gcdd}
\end{figure}

The involution depends on the following factorisation of the paths in $R_2^2$.  
%===== Lemmma
\begin{lemma} \label{lem:34fact}
Let $\om\in  R_2^2$,      then $\om$  can be factorised in one and only one   of the five  following forms (illustrated in \figref{fig:gcdd}):
\begin{subequations}\label{eq:f5}
\begin{align} 
\om^{(1)}& = A_1 B_1\, d u M_1 \, d'\,  L_{1}		&  w(d)   =\bb, \quad w(d')=\ab	\label{eq:f1}\\
\om^{(2)} &=  A_2 B_2 \, u M_2 d d'\,  L_{2}\,		&  w(d')  =-\ab\bb				\label{eq:f2}\\
\om^{(3)}&= A_3 B_3\,  u\, d\,  L_{3},				&  w(d)  = -\ab\bb				\label{eq:f3}\\
\om^{(4)}& =\, A_4 B_4 d u\,  L_{4} ,				&  w(u)  =\ab\bb 				\label{eq:f4}\\
\om^{(5)}& =\, A_5 B_5  & 															\label{eq:fp}
\end{align}
where $u$, $u'$ are up steps, $d$, $d'$ are down steps, $A_i$, $B_i d$, $M_i$ and $u L_id$ are all (possibly empty) elevated Dyck subpaths and $w(e_n)$ is the weight of step $e_n$. The subpaths $A_i$   contain  only $\ab$ weighted steps, the subpaths $B_i$  contain only $\bb$ weighted steps, the subpaths $M_i$   contain  no weighted steps and the subpaths $L_i$   contain any weighted steps (ie.\ $\ab$, $\bb$, $\ab\bb$ and $-\ab\bb$).
\end{subequations}
\end{lemma} 
%=====
The factorisation is defined by what will be referred to as a ``\textbf{bad}'' step. Bad steps (if they occur) are of two types:  1) an `$\ab\bb$-bad' step or 2) an `$\ab$-bad' step. An $\ab\bb$-bad step is the leftmost step weighted $\pm\ab\bb$ and an $\ab$-bad step is the   leftmost  step weighted $\ab$ occurring to the right of a step weighted $\bb$. Note, the $R_4$ paths are precisely the paths with no bad steps. The factorisation  cases are as follows:
\begin{itemize}

\item The path has a bad step:
\begin{itemize}
\item The leftmost bad step is an $\ab$-step. Thus to the left of the $\ab$-step there are no $\ab\bb$ weighted steps and hence the  path must factor according to case \eqref{eq:f1}.

\item The leftmost bad step is an $\ab\bb$-step. There are two sub-cases:
\begin{itemize}
\item The $\ab\bb$-step is above height one (ie.\ negative). We split this into two further cases depending on:
\begin{itemize}
\item whether the step before the bad step is a  down step -- case \eqref{eq:f2}
\item or an up step -- case \eqref{eq:f3}. 
\end{itemize}
\item The $\ab\bb$-step is below height one (ie.\ positive). This is case \eqref{eq:f4}.
\end{itemize}

\end{itemize}
\item The path has no bad step -- thus contains no $\ab\bb$ steps and all the $\ab$ steps are to the left of the $\bb$ steps. This is case \eqref{eq:fp}.
\end{itemize}

The involution $\Phi_3 $,  detailed below, can be succinctly summarised as follows.
Referring to \figref{fig:gcdd}: In (a) flip the pair of edges to the left of $M_1$,  one of which is now a down edge and move this one to the other side of $M_1$ together with  the factor $\bb$ (and change its sign). This is now the same as (b).  In (c)  flip the pair of edges to the left of $M_1$ and change the sign giving and (d). Hence (a) and (b) cancel as do  (c) and (d) leaving only (e).

The involution $\Phi_3 $  is defined on the path set $R_3^2$ and will have fixed point set $R_4$. Define the signed set as follows. Let
\begin{align}
	R_3^2&=\Om^{(3)}_{+}\cup \Om^{(3)}_{-}\\
\intertext{where the signed sets are}
	\Om^{(3)}_{+}&=\{\om\suchthat  \text{$\om\in R_2^2$ and $\om$ has positive weight} \}\\
	\Om^{(3)}_{-}&=\{\om\suchthat  \text{$\om\in R_2^2$ and and $\om$ has negative weight}\}.
\end{align}
The involution $\Phi_3: R_3^2\to R_3^2$,   falls into five cases corresponding to the five factorisations. Let $\om\in R_3^2$ and  $\om'=\Phi_3(\om)$.

\begin{enumerate}

% Case 1 =====
\item  If $\om$ is of the  form of  \eqref{eq:f1}  then $\om'$ is obtained from $\om$ by moving   $d$ to the right of $d'$, removing the $\ab$ and $\bb$ weights from $d$ and $d'$ and giving the moved  $d$ step weight $-\ab\bb$, that is,
\begin{equation}  	
		A_1 B_1\, d u M_1 \, d'\,  L_{1}	 \longrightarrow \om'= A_1B_1 \, u M_1 d'  d \,  L_{1} 
\end{equation}
Thus $\om' $ is of the form and weight of \eqref{eq:f2}. In $\om'$, $w(d)w(d') =-\ab\bb$ and thus the sign of $\om'$ is opposite to that of $\om$ as required.

% Case 2 =====
\item If $\om$ is of the form of \eqref{eq:f2} then $\om'$ is obtained from $\om$ by shifting $d$ to the left of $u$, changing the weight of $d$ to $\bb$, and that of $d' $ to $\ab$, to give
\begin{equation}   
A_2 B_2 \, u M_2 d d'\,  L_{2}\to \om'  = A_2 B_2\, d  u  M_2 \, d'\, L_{2}
\end{equation}
Thus $\om' $ is of the form and weight of  \eqref{eq:f1}.  Since now $w(d)w(d') =+\ab\bb$ the sign of $\om'$ is opposite to that of $\om$ as required.

% Case 3 =====
\item If $\om$ is of the form of \eqref{eq:f3} then $\om'$ is obtained from $\om$ by swapping the $u$ and $d$ steps and changing the weight of $d $ to $+\ab\bb$, to give
\begin{equation} 
	A_3 B_3\,  u\, d\,  L_{3}\longrightarrow \om'=A_3 B_3\, d\,  u \,  L_{3}
\end{equation}
Thus $\om' $ is of  the form \eqref{eq:f4}. Since now $w(  d )=+\ab\bb$ the sign of $\om'$ is opposite to that of $\om$ as required.

% Case 4 ====
\item  If $\om$ is of the form of \eqref{eq:f4} then $\om'$ is obtained from $\om$ by swapping the $u$ and $d$ steps and changing the weight of $d $ to $-\ab\bb$ to give
\begin{equation}
		A_4 B_4\, d\, u\,  L_{4}\longrightarrow\om'   = A_4 B_4\,  u \, d \,  L_{4}
\end{equation}
Thus $\om' $ is of  the form \eqref{eq:f3}. Since now $w(  d )=-\ab\bb$ the sign of $\om'$ is opposite to that of $\om$ as required.

% Case 5 ====
\item If $\om$ is of the  form of \eqref{eq:f5} then $\om'=\om$. This is the fixed point set.

\end{enumerate}

In all cases after the action of $\Phi_3$, the bad step stays immediately to the left of the initial $L_i$ factor thus ensuring $\Phi_3^2=1$ as required.  The fixed point set has no bad steps ie.\ all the $\ab$ weighted steps are to the left of the $\bb$ steps and there are no $\pm\ab\bb$ weighted steps -- thus the fixed point set is the set $R_4$ as desired.

\newpage

%================================================================================
 \bibliographystyle{plain} 
 %================================================================================

\bibliography{bibliography} % the path to the .bib file is kept in the macros.tex file

\end{document}